\documentclass[a4paper]{article}

\usepackage{cmap}
\usepackage[T2A]{fontenc}
\usepackage[utf8]{inputenc} % [utf8][cp866][cp1251]
\usepackage[english]{babel}
\usepackage{amsmath, amssymb, latexsym}
\usepackage{tikz}
\usepackage{ifthen}
\usepackage{hyperref}
\usepackage{longtable}
\usepackage{setspace}
\usepackage{indentfirst}
\title{On multifold perfect codes and some other completely regular codes in the Doob graphs and quaternary Hamming graphs
\thanks{The work was funded by the Russian Science Foundation, Grant 22-11-00266}}
\author{Evgeny Bespalov
\thanks{Sobolev Institute of Mathematics, Novosibirsk, Russia. E-mail: bespalovpes@mail.ru}}

\date{}

\makeatletter
\def\@seccntformat#1{\csname the#1\endcsname.\ } % точка после номера раздела
\makeatother

\newif\ifNoRemark
\def\addtheorem#1#2#3#4{
\ifthenelse{\equal{#2}{}}{}%
{\ifthenelse{\expandafter\isundefined\csname the#2\endcsname}{\newcounter{#2}}{}}
\newenvironment{#1}[1][\global\NoRemarktrue]% No Remark by default
{\par\addvspace{2mm plus 0.5mm minus 0.2mm}\noindent % new paragraph without indent
{\bf #3}\ifthenelse{\equal{#2}{}}{}%
{\refstepcounter{#2}{\bf ~\csname the#2\endcsname}}%
{\bf \vphantom{##1}\ifNoRemark.\ \else\ (##1).\fi}\begingroup #4}%
{\endgroup\par\addvspace{1mm plus 0.5mm minus 0.2mm}\global\NoRemarkfalse}
\expandafter\newcommand\csname b#1\endcsname{\begin{#1}}
\expandafter\newcommand\csname e#1\endcsname{\end{#1}}
}

\addtheorem{theoreman}{thrm}{Theorem}{\sl}
\addtheorem{lemman}{lmm}{Lemma}{\sl}
\addtheorem{predln}{prpstn}{Proposition}{\sl}
\addtheorem{coroll}{crllr}{Corollary}{\sl}
\addtheorem{remarkn}{rmrk}{Remark}{}
\addtheorem{definitionn}{def}{Definition}{}

 \newenvironment{proof}[1][\hspace{-1.0ex}]%
  {\par\addvspace{1mm}{\sc Proof\hspace{1.0ex}{#1}.} }%
  {\quad$\blacktriangle$\par\addvspace{1mm}}
	
\def\shpart#1 #2 #3 #4!{node [#1] {} +(1,0) node [#2] {} +(2,0) node [#3] {} +(3,0) node [#4] {}}
\def\sh#1#2#3#4{\begin{tikzpicture}[
scale=0.7,
nz/.style={circle,fill=white,draw=black, 
           inner sep=2.5pt},
xz/.style={circle,fill=black!50!white,draw=black, 
           inner sep=2.5pt},
vz/.style={circle,fill=black!23!white,draw=black, 
           inner sep=2.5pt},
wz/.style={circle,fill=black!66!white,draw=black, 
           inner sep=2.5pt},
zz/.style={circle,fill=black,draw=black, 
           inner sep=2.5pt},
yz/.style={circle,fill=yellow,draw=black, 
           inner sep=2.5pt},
rz/.style={circle,fill=red,draw=black, 
           inner sep=2.5pt},
bz/.style={circle,fill=blue,draw=black, 
           inner sep=2.5pt},
gz/.style={circle,fill=green,draw=black, 
           inner sep=2.5pt},  
oz/.style={circle,fill=orange,draw=black, 
           inner sep=2.5pt},
fz/.style={circle,fill=violet,draw=black, 
           inner sep=2.5pt}, 
kz/.style={circle,fill=brown,draw=black, 
           inner sep=2.5pt},  
cz/.style={circle,fill=cyan,draw=black, 
           inner sep=2.5pt},           
scale=0.7]
\begin{scope}
\clip [xslant=-0.577] (-1.4,-1.20) rectangle (2.4,2.1);
\draw[xslant=0.577,ystep=.866,xstep=1,draw=black] (-4.9,-2.1) grid (5.4,3.9);
\draw[xslant=-0.577,ystep=9.866,xstep=1,draw=black] (-3.4,-2.1) grid (6.4,3.9);
\draw (-120:1) \shpart #4!
++(120:1) \shpart #3!
++(120:1) \shpart #2!
++(120:1) \shpart #1!;
\end{scope}
\end{tikzpicture}}

\def\shpart#1 #2 #3 #4!{node [#1] {} +(1,0) node [#2] {} +(2,0) node [#3] {} +(3,0) node [#4] {}}
\def\hh#1#2#3#4{\begin{tikzpicture}[
scale=0.7,
nz/.style={circle,fill=white,draw=black, 
           inner sep=2.5pt},
xz/.style={circle,fill=black!50!white,draw=black, 
           inner sep=2.5pt},
vz/.style={circle,fill=black!23!white,draw=black, 
           inner sep=2.5pt},
wz/.style={circle,fill=black!66!white,draw=black, 
           inner sep=2.5pt},
zz/.style={circle,fill=black,draw=black, 
           inner sep=2.5pt},
yz/.style={circle,fill=yellow,draw=black, 
           inner sep=2.5pt},
rz/.style={circle,fill=red,draw=black, 
           inner sep=2.5pt},
bz/.style={circle,fill=blue,draw=black, 
           inner sep=2.5pt},
gz/.style={circle,fill=green,draw=black, 
           inner sep=2.5pt},  
oz/.style={circle,fill=orange,draw=black, 
           inner sep=2.5pt},
fz/.style={circle,fill=violet,draw=black, 
           inner sep=2.5pt}, 
kz/.style={circle,fill=brown,draw=black, 
           inner sep=2.5pt},  
cz/.style={circle,fill=cyan,draw=black, 
           inner sep=2.5pt},           
scale=0.7]
\begin{scope}
\clip [xslant=-0.577] (-1.4,-1.20) rectangle (2.4,2.1);
\draw[xslant=0.577,ystep=.866,xstep=1,draw=black] (-4.9,-2.1) rectangle (5.4,3.9);
\draw[xslant=-0.577,ystep=9.866,xstep=1,draw=black] (-3.4,-2.1) rectangle (6.4,3.9);
\draw (-120:1) \shpart #4!
++(120:1) \shpart #3!
++(120:1) \shpart #2!
++(120:1) \shpart #1!;
\end{scope}
\end{tikzpicture}}

\begin{document}

\maketitle

\begin{abstract}
We consider the problem of existence of perfect $2$-colorings in the Doob graphs $D(m,n)$ and $4$-ary Hamming graphs $H(n,4)$. We characterize all parameters for which multifold $1$-perfect code in $D(m,n)$ exists. 
Also, we prove that for any pair $(b,c)$ that satisfy standard conditions (Lloyd's and sphere-packing conditions) there is perfect $(b,c)$-coloring in Doob graphs and Hamming graphs if diameter of graph and $n$ are sufficiently large ($n$ not less than $0$, $1$ or $8$ for some cases).  Also we obtain some completely regular codes with covering radius $2$.
\end{abstract}

\def\VV{{\scriptscriptstyle\mathrm{V}}}

\section{Introduction}

An equitable partition of graph is a partition of its vertex set into $k$ cells such that induced subgraph on every cell is a regular graph and bipartite graph generated by edges between any two cells is biregular. In the current paper we use equivalent object, a perfect $k$-coloring, a $k$-valued function on a graph such that preimaiges of values form equitable partitions. 
A code is completely regular with covering radius $\rho$ if its distance coloring is a perfect $(\rho+1)$-coloring. In case of two colors a set of first (second) color in perfect $2$-coloring is a completely regular code with covering radius $1$. 

Many combinatorial objects are correspond to perfect colorings, for example $1$-perfect codes, latin squares and hypercubes. In some cases, objects that attain some bound, for example orthogonal arrays that lie on Bierbrauer-Friedman bound, are equivalent to perfect colorings. Some more example are given in \cite{PotAvg:perfcol}. So, there are many work with study of perfect colorings in many graphs. For perfect $2$-colorings in Hamming graphs $H(n,q)$ we refer to \cite{FDF:PerfCol} (binary case),  \cite{BKMTV:perfcolinham} (case when $q$ is a prime power) and \cite{BRZ:crg} (overview on completely regular codes). Also perfect colorings studied in Johnson graphs \cite{EGGV:johnsonperf}, \cite{Mog:crgjohn}, non-distance regular graphs \cite{Puz:2011en} and other graphs.

One of the most important completely regular codes are $1$-perfect codes.  
The code is called $1$-perfect if any radius-$1$ ball contains $1$ vertex from the code. If any radius $1$-ball has $\mu$ vertices from the code then such code is called $\mu$-fold $1$-perfect code.
%is generalization of $1$-perfect code. 
The problem of existence of multifold $1$-perfect codes in  Hamming graphs on prime power alphabet was completely solved in \cite{Kro:multperf}. It was proved that such codes exist if and only if Lloyd's condition and sphere-packing condition are hold. In the current paper, we obtain the same result for the Doob graphs.

Besides multifold $1$-perfect codes we consider the problem of the existence of other perfect $2$-colorings (we also denote them as perfect $(b,c)$-colorings where $b$ and $c$ are degrees in bipartite biregular graph obtained by edges between first and second colors) 
in quaternary Hamming graphs $H(N,4)$ and in Doob graphs $D(m,n)$ which are distance-regular with the same parameters as Hamming graph $H(2m+n,4)$. 
Some codes that correspond to perfect $2$-colorings have already been studied in Doob graphs. The study of $1$-perfect codes started in \cite{KoolMun:2000}. In \cite{Kro:pfdoob} all $m,n$ for which there is a $1$-perfect code in $D(m,n)$ were completely characterized. In \cite{KroBes:MDS2} all $2$-MDS codes were characterized.

Recall (see \cite{BKMTV:perfcolinham}) that if there is a perfect $(b,c)$-coloring in Hamming graph $H(N,q)$ then perfect $(b,c)$-coloring also exist in $H(n,q)$ for any $n \ge N$. So, problem of existence of perfect $2$-colorings can be divided into two problems: to characterize all admissible pairs $(b,c)$ (i.e. such that there is perfect $(b,c)$-coloring in $H(N,q)$ for some $N$) and for any admissible pair $(b,c)$ to find minimum such $N$. 
For Doob graphs the definition of admissibility is slightly more difficult because for fixed diameter there are a few Doob graphs of this diameter. So we define admissibility for Doob graphs that depends on diameter with restrictions on $n$. See more details in Subsection~\ref{s:nec}.

In this paper we generalized the ideas from \cite{FDF:PerfCol} and \cite{BKMTV:perfcolinham} that used in constructions a partition of $2$-MDS code in $H(q,q)$ into codes with distance $3$. We construct such partitions for $2$-MDS codes in some other quaternary Hamming graphs and Doob graphs. We should note that such partitions for Hamming graphs can be also obtained using a cosets of additive completely regular codes that dual to known additive two-weight codes from \cite{SZZ:equidistance}, \cite{BZZ:cover2} (see Remark~\ref{r:mdspart} for more details). In the result, we complete characterization of all admissible pairs $(b,c)$ for $4$-ary Hamming graphs considering the remain case ($b$ and $c$ are relatively prime and $(b+c)/gcd(b,c)$ is a odd power of $2$) after works $\cite{BKMTV:perfcolinham}$, \cite{Kro:multperf}. For Doob graphs we obtain the same result except case when $n$ is small for some pairs $(b,c)$. For more details see Theorem~\ref{t:admissible}. Also we construct some completely regular codes with covering radius $2$ that obtained from the same code in $4$-ary Hamming graphs (see Corollary~\ref{c:rad2}).

\section{Preliminaries}

\subsection{Basic definitions}

Given a graph $G$ denote by $\VV{G}$ its vertex set. A \emph{code} $C$ in $G$ is an arbitrary nonempty set of $\VV{G}$. If graph is connected for vertices $x$ and $y$ we define \emph{distance $d(x,y)$} between $x$ and $y$ as a length of the shortest path between them. A \emph{code distance} of $C$ is a minimum distance between two different vertices from $C$, i.e. $\min_{x \ne y, x,y \in C}\{d(x,y)\}$. The covering radius $\rho(C)$ is a maximum possible distance between $C$ and vertex of $G$, i.e. $\rho(C)=max_{x \in \VV{G}}(d(x,C))$.  

A surjective function $f$ of a graph $G$ into color set of cardinality $k$ is called \emph{$k$-coloring}, the elements often denoted by $\{1,\ldots,k\}$ are called \emph{colors}. If for any colors $i$ and $j$ (not necessary different) an arbitrary vertex $x$ of color $i$ adjacent to exactly $s_{i,j}$ vertices of color $j$, where $s_{i,j}$ does not depend on choice of $x$, the coloring $f$ is called \emph{perfect}. 
The matrix $S=(s_{i,j})$ of order $k$ is called the \emph{quotient matrix} of $f$. 
The correspondence partition of $\VV{G}$, i.e. $(f^{-1}(1),\ldots,f^{-1}(k))$, is called \emph{equitable partition} with quotient matrix $S$. 

A code $C$ with covering radius $\rho$ is called \emph{completely regular}  if a distance $(\rho(C)+1)$-coloring $f$, i.e. $f(x)=i$ if $d(x,C)=i$, is a perfect coloring with tridiagonal quotient matrix. The array 
$[s_{0,1},\ldots,s_{\rho-1,\rho},s_{1,0},\ldots,s_{\rho,\rho-1}]$ is called \emph{intersection array} of completely regular code $C$. 
Graph $G$ is called \emph{distance-regular} if any vertex of $G$ is a completely regular code with same quotient matrix. In case of $2$ colors, for any color $i \in \{1,2\}$ a color set $f^{-1}(i)$ of perfect $2$-coloring is a completely regular code with covering radius $1$. A perfect $2$-coloring with quotient matrix $S$ we also will denote as perfect $(s_{1,2},s_{2,1})$-coloring.        

The Hamming graph $H(n,q)$ is a direct product of $n$ copies of complete graph $K_q$ on $q$ vertices. 
The Doob graph $D(m,n)$ is a direct product of $m$ copies of the Shrikhande graph and $n$ copies of complete graph $K_4$ on $4$ vertices. The Shrikhande graph can be represented as Cayley graph on group $\mathbb{Z}^2_4$ with connecting set $\{01,03,10,30,11,33\}$ (the vertices of graph are elements of $\mathbb{Z}^2_4$ and two vertices are adjacent if and only if their difference belongs to connecting set). The vertex set of $D(m,n)$ can be represented as $(\mathbb{Z}^2_4)^m \times (\mathbb{Z}_4)^n=
\{(x^{*}_1,\ldots,x^{*}_m;x'_1,\ldots,x'_n: x^{*}_i \in \mathbb{Z}^2_4, x'_j \in \mathbb{Z}_4)\}$. For any vertex $x$ of $D(m,n)$ we denote it by $(x^{*}_1,\ldots,x^{*}_m;x'_1,\ldots,x'_n)$. where $x^{*}_i \in \mathbb{Z}^2_4$ and $x'_j \in \mathbb{Z}_4$.
If $m>0$ graph $D(m,n)$ is called \emph{Doob graph}, if $m=0$ the graph $D(0,n)$ is the Hamming graph $H(n,4)$. 

The graph $D(m,n)$ is a distance regular with the same parameters as a Hamming graph $H(2m+n,4)$. The graph $D(m,n)$ has the following eigenvalues: $6m+3n-4i, i=0,1,\ldots,2m+n$.

Denote by $B(x)$ a radius-$1$ ball in $D(m,n)$ with center $x$, 
i.e. $B(x)=\{y \in \VV{D(m,n)}: d(x,y) \le 1\}$. A code $C$ is called \emph{$\mu$-fold $1$-perfect} if for any vertex $x$ of $D(m,n)$ there are $\mu$ vertices of $C$ in $B(x)$. If $\mu=1$, $C$ is called \emph{$1$-perfect code}. Note that $C$ is $\mu$-fold $1$-perfect code in $D(m,n)$ if and only if $C$ is a completely regular code with intersection array $[6m+3n-\mu+1,\mu]$.  
A code is called \emph{$2$-MDS code} if $C$ is an independent set of cardinality $4^{2m+n-1}$. A code $C$ is a $2$-MDS code in $D(m,n)$ if and only if $C$ is completely regular with intersection array $[6m+3n,2m+n]$.

Denote by $J_{k \times l}$ the $k \times l$-matrix such that any element equals $1$. Denote by $E_{t \times t}$ the identity $t \times t$-matrix. We will omit indices if it is clear from context.

\subsection{Necessary conditions}\label{s:nec}

\begin{predln}\label{p:admiss}
Let there is a perfect coloring $g$ in $D(m,n)$ with quotient matrix $S$. 
Then there is a perfect coloring $f$ in $D(m+m',n+n')$ with quotient matrix $S+(6m'+3n')E$ for any non-negative integer $m'$ and $n'$.
\end{predln}
\begin{proof}
Define the function
$f:\VV{D(m+m',n+n')} \to \VV{D(m,n)}$ in the following way:
$$f(x^{*}_1,\ldots,x^{*}_{m+m'};x'_1,\ldots,x'_{n+n'})=g(x^{*}_1,\ldots,x^{*}_m;x'_1,\ldots,x'_n).$$

Let $x$ be an arbitrary vertex of $D(m+m',n+n')$ that has color $i$. Represent the neighbourhood $N(x)$ of $x$ as $N_1(x) \cup N_2(x)$, where $N_1(x)$ consists of neighbours that differs from $x$ in some position 
from $\{1,\ldots,m,m+m'+1,\ldots,m+m'+n\}$, and $N_2(x)$ consists from neighbours that differs in some other position. By the definition, $N_1(x)$ have $s_{i,j}$ vertices of color $j$ and all vertices of $N_2(x)$ has the same color as $x$. So, $x$ has $s_{i,j}$ neighbours of color $j$ if $i \ne j$ and $(s_{i,i}+6m'+3n')$ if $i=j$.
\end{proof}

Proposition~\ref{p:admiss} allows us to define admissible arrays.
Let $a$ be a non-negative integer. Intersection array $[b;c]$ is called 
\emph{$a$-admissible} if there is $D$ such that there is $(b,c)$-coloring in $D(m,n)$ for any $m,n$ where $2m+n \ge D$, $n \ge a$.
Intersection array $[b,c]$ is called 
\emph{$\infty$-admissible} if there is $N$ such that there is $(b,c)$-coloring in $H(n,4)$ for any $n$ where $n \ge N$.

\begin{predln}\label{p:necc1}
If there is $(b,1)$-coloring in $D(m,n)$, then $b+1$ divides by $3$. Moreover, if $n=0$, then $b$ divides by $6$.
\end{predln}
\begin{proof}
Consider an arbitrary vertex $x$ of color $1$. The neighbourhood of $x$ is a union of $m$ $6$-cycles and $n$ $3$-cycles. The vertices in any cycle has the same color. Indeed, if there two vertices of different colors in some cycle then there are some adjacent vertices $u$ and $v$, such that $u$ has color $1$ and $v$ has color $2$. But in this case, $v$ has two neighbours of color $1$ and we have a contradiction with $c=1$. Hence, any $6$-cycle contains $0$ or $6$ vertices of color $2$ and any $3$ cycle contains $0$ or $3$ vertices of color $2$. So, $b$ is divisible by $3$, moreover, if $n=0$ then $b$ is divisible by $6$.
\end{proof}

\begin{predln}\label{p:neceigen}
    If there exists a $(b,c)$-coloring in $D(m,n)$, then 
    \begin{enumerate}
    \item $\displaystyle{\frac{b+c}{gcd(b,c)}=2^k}$ for some positive $k$.
    \item $b+c=4i$ for some $i \in \{1,\ldots,2m+n\}$.
    \item if $c=1$ ($b=1$) then $b=4^l-1$ ($c=4^l-1$) for some $l$, moreover, $n\ne 0$.
    \end{enumerate}
\end{predln}
\begin{proof}
1) By counting the edges between vertices of different colors we have 
$|f^{-1}(1)|b=|f^{-1}(2)|c$, and hence, $f^{-1}(1)=\displaystyle{\frac{c}{b+c} \cdot 4^{2m+n}}$. Since the cardinality of $f^{-1}(1)$ is integer, $\frac{b+c}{gcd(b,c)}$ is a power of $2$.

2) The quotient matrix of $f$ has eigenvalue $6m+3n$ and $6m+3n-(b+c)$. By Lloyd's condition these values should be eigenvalues of the graph, which follows $(b+c)=4i$ for some $i \in \{1,\ldots,2m+n\}$.

3) By i.1 we have $b+1=2^k$ for some integer $k$. By Proposition~\ref{p:necc1} we have that $b$ is divisible by $3$ and if $n=0$ then $b$ is divisible by $6$.  Hence, $k$ is even and $n \ne 0$. 
\end{proof}

\subsection{1-perfect codes and 2-MDS codes}

In \cite{Kro:pfdoob} the problem of existence of $1$-perfect codes in the Doob graphs was completely solved.  

\begin{theoreman}\cite{Kro:pfdoob}\label{t:perfcodeexist}
There is a $1$-perfect code in $D(m,n)$ if and only if $2m+n=\displaystyle{\frac{4^{l}-1}{3}}$ for some positive integer $l$.
\end{theoreman}

A partition of $D(m,n)$ into $1$-perfect codes give the following perfect coloring.

\begin{coroll}\label{c:JE}
If $2m+n=\displaystyle{\frac{4^l-1}{3}}$ then there is a perfect $4^l$-coloring with quotient matrix $(J-E)$.
\end{coroll}
\begin{proof}
This perfect coloring can be obtained using a partition into $4^l$ disjoint $1$-perfect codes.
\end{proof}

A partitions of $D(m,n)$ into disjoint $2$-MDS codes gives the following perfect colorings. 

\begin{predln}\label{p:latin}
For any  $m, n$ there is a perfect $4$-coloring in $D(m,n)$ with quotient matrix $(2m+n)(J-E)$.
\end{predln}

A perfect $2^{k+2}$-coloring is called \emph{$k$-multipartite perfect coloring} if its colors can be represented as $(i,j)$, $i \in \{0,1,2,3\}$, 
$j \in \{0,1,\ldots,2^k-1\}$, such that for any $i,j$ any vertex of color $(i,j)$ has exactly $1$ neighbour of color $(r,s)$ if $(r \ne i)$ and $0$ neighbours of color $(r,s)$ if $r=i$. Note that by condition all color sets have the same cardinality, so the number of colors should be a power of $2$. Also note that the union of colors $(i,0),\ldots,(i,2^k-1)$ is a $2$-MDS code for any $i$. 

\begin{lemman}\label{l:multipart}
Let there is a partition of $\VV{D(m,n)}$, where $2m+n=2^k$ for some positive integer $k$, into disjoint codes $C^0_0,\ldots,C^0_{2^k-1}, C^1_{0},\ldots,C^3_{2^k-1}$ such that $C^i_j$ has code distance $3$, $C^i_0 \cup \ldots \cup C^i_{2^k-1}$  is a $2$-MDS code for any $i,j$. Then coloring $f$ with colors $(0,0),\ldots,(3,2^k-1)$ such that $f(x)=(i,j)$ if $x \in C^i_j$ is $k$-multipartite perfect coloring in $D(m,n)$. 
\end{lemman}
\begin{proof}
%Denote colors of $2^{k+2}$ coloring as $(i,j)$, where $i \in \{0,\ldots,3\}$, $j \in \{0,\ldots,2^k\}$.
%Define coloring $f$ as follows: $f(x)=(i,j)$ if $x \in C^i_j$. 
Let $x$ be an arbitrary vertex of color $(i,j)$. Let us count the number of neighbours of color $(r,s)$. If $r=i$ then this equals $0$ because 
$C^i_j$ and $C^i_s$ belong to the same $2$-MDS code.
Let $i \ne r$. The vertex $x$ adjacent to $2^k$ vertices of 
$C^r_0 \cup \ldots \cup C^r_{2^k}$ since this union is $2$-MDS code. Since $C^r_j$ for any $j$ has a code distance at least $3$ the vertex $x$ can't have more than $1$ neighbour in $C^r_l$ for any $l$. So, $x$ has exactly $1$ neighbour of color $(r,s)$ and $f$ is $k$-multipartite perfect coloring. 
\end{proof}

%\begin{lemman}\label{l:mult3j}
%Let there is a $k$-multipart perfect coloring in $D(m,n)$, where $2m+n=2^k$. %Then there is a perfect $2^k$-coloring in $D(m,n)$ with quotient matrix %$3J$. Moreover, if $n>0$ then there is a perfect $2^k$-coloring in $D(m,n-%1)$ with quotient matrix $3(J-E)$.
%\end{lemman}
%\begin{proof}
%
%\end{proof}

\begin{predln}\label{p:multdiam4}
There is a $2$-multipartite perfect coloring in $D(m,n)$, where $2m+n=4$. 
\end{predln}
\begin{proof}
By Theorem~\ref{t:perfcodeexist} there is a $1$-perfect code $C$ in $D(m,n+1)$. 
Denote $M^0_j$ the set of vertices with $j$ in the last position. 
$M^0_j$ has code distance $3$ for any $j$ and $M^0_0 \cup \ldots \cup M^0_3$ is $2$-MDS code. If $n=2$ denote $M^i_j=\{x+(00,0,i):x \in M^0_j\}$. The codes $M^0_0,\ldots,M^3_3$ satisfy the condition of Lemma~\ref{l:multipart} so the statement is true. The case $n=4$ is analogous. Let $n=0$. Denote $M^1_j=\{x+(00,01):x \in M^0_j\}$, $M^2_j=\{x+(00,10):x \in M^0_j\}$, $M^3_j=\{x+(00,11):x \in M^0_j\}$. Note that there are $2$ different $1$-perfect codes in $D(2,1)$ up to equivalence and they are explicitly described (see \cite{BesKro:mdsdoob}, Appendix). In both cases in our construction $2$-MDS codes $M^1_0 \cup \ldots \cup M^1_3$ and $M^3_0 \cup \ldots \cup M^3_3$ are disjoint. So, the codes $M^0_0,\ldots,M^3_3$ satisfy the condition of Lemma~\ref{l:multipart} so the statement is true.
\end{proof}

\begin{lemman}\label{l:3jmult}
Let there is a $k$-multipartite perfect coloring in $D(m,n)$. Then there is a perfect perfect $2^k$-coloring with quotient matrix $3J$
\end{lemman}
\begin{proof}
Let $g$ be a $k$-multipartite perfect coloring in colors $(0,0),
\ldots,(3,2^k-1)$.
The requirement coloring can be obtained if for any $j$ we unit colors $(0,j)$, $(1,j)$, $(2,j)$ and $(3,j)$. 
\end{proof}

\section{Colorings of direct products}

The $n \times n$ matrix $M=(m_{i,j})$ is called \emph{equal-diagonal} if $m_{i,j}=m_{i+1,j+1}$, where addition is modulo $n$, for any $i,j$.

\begin{theoreman}\label{t:diag}
    Let $g$ be a perfect $k$-coloring of graph $G$ with quotient matrix  \[S=\displaystyle{
    \begin{pmatrix}
A_{1,1} & A_{1,2} & \ldots & A_{1,l} \\   
A_{2,1} & A_{2,2} & \ldots & A_{2,l} \\
\ldots & \ldots & \ldots & \ldots \\
A_{l,1} & A_{l,2} & \ldots & A_{l,l} 
    \end{pmatrix}
    }\] 
    where $A_{i,i}$ is equals-diagonal $|k_i| \times |k_i|$ for any $i$,  $A_{i,j}$ is the matrix $a_{i,j}J_{|k_i| \times |k_j|}$ and $k_1+\ldots+k_l=k$

    Let $h_i$ be a perfect $k_i$-coloring of graph $H$ with quotient matrix $B_i$ that is equal-diagonal.

    Then there is a perfect coloring $f$ of graph $F=G \times H$  with quotient matrix \[\displaystyle{
    \begin{pmatrix}
A_{1,1}+B_1 & A_{1,2} & \ldots & A_{1,l} \\   
A_{2,1} & A_{2,2}+B_2 & \ldots & A_{2,l} \\
\ldots & \ldots & \ldots & \ldots \\
A_{l,1} & A_{l,2} & \ldots & A_{l,l}+B_l 
    \end{pmatrix}
    }\]  
\end{theoreman}
\begin{proof}
   Represent the colors of $k$-coloring as $(i,j)$ where $i \in \{1,\ldots,l\}$, $j \in \{1,\ldots,k_i\}$, $k_1+\ldots+k_l=k$.
    Denote the vertex of $F$ as $z=(x,y)$, where $x \in \VV{G}$, $y \in \VV{H}$.
    Define the coloring $f$ in the following way:
    $$f(x,y)=(i,h_i(y)+j \mod k_i)$$
    where $i,j$ such that $g(x)=(i,j)$.

Let $z=(x,y)$ be an arbitrary vertex of $F$ of color $(i,s)$, where   
    $f(x)=(i,j)$ and $h_i(y)+j=s$.
Let us count the number of vertices of color $(r,m)$ that adjacent to $z$.     

Let $i \ne r$. Then $z$ has $b_{i,r}$ neighbours of color $(r,m)$ for any $m \in \{1,\ldots,k_r\}$. 
Indeed, a vertex $z'=(x',y')$, where $g(x')=(r,t)$, of color $(r,m)$ adjacent to $z$ if and only if $x$ and $x'$ are adjacent in $G$ and $t=m-g_r(y')$. So, $z$ has number of neighbours of color $(r,m)$ equal to the number of neighbours of color $(r,t)$ of $x$ in graph $G$. This number equals $a_{i,r}$. 

Let $i=r$. 
Represent the neighbourhood $N(x)$ of $x$ as $N_1(x) \cup N_2(x)$, where $N_1(x)$ is a set of neighbours $(x,y')$ such that $y$ and $y'$ are adjacent in $H$ and $N_2(x)$ is a set of neighbours $(x',y)$ such that $x$ and $x'$ are adjacent in $G$. A vertex $z_1=(x,y')$, where $y'$ adjacent to $y$, has color $(i,m)$ if and only if $m=h_i(y')+j$. Hence, $|N_1(x)|$ equals $(s-j,m-j)$-element of matrix $B_i$. Since $B_i$ is equal-distance we have that $|N_1(x)|$ equals $(s,m)$-element of $B_i$. A vertex $z_2=(x',y)$, where $g(x')=(i,t)$, has color $(i,m)$ if and only if $h_i(y)+t=m$. Hence, $|N_2(x)|$ equals $(s-h_i(y),m-h_i(y))$-element of matrix $A_{i,i}$. Since $A_{i,i}$ is equal-distance we have that $|N_2(x)|$ equals $(s,m)$-element of $A_{i,i}$ 
%A vertex $z'=(x',y')$ of color $(i,m)$, where $g(x')=(i,t)$, adjacent to $z$ if and only if either $x=x'$, $y$ and $y'$ are adjacent in $H$ and $h_i(y')+j=m$ or $y=y'$, $x$ and $x'$ are adjacent in $G$ and $f(x')=(i,j+m-s)$.   
%Since $g_i(y)+j=s$ and $f(x)=(i,j)$ we have that the number of vertices of color $(r,m)$ that adjacent to $z$ equals the number of $b^i_{s-j,m-j}+a^i_{j,m-s+j}$, where $a^i_{u,v}$ and $b^i_{u,v}$ are elements of $A_{i,i}$ and $B_i$ correspondingly. Since $A_{ii}$ and $B_{i}$ are equal-diagonal this equals $b^i_{s,m}+a^i_{s,m}$. This complete the proof.
\end{proof}

\begin{predln}\label{p:multiplyk}
Let there is a perfect coloring $g$ in $D(m,n)$ with quotient matrix $S$. Then
\begin{enumerate}
    \item  For any even $k$ there is a perfect coloring in $D(mk+m',n')$, where $2m'+n'=kn$, with quotient matrix $kS$.
    \item  For any odd $k$ there is a perfect coloring in $D(mk+m',n')$, where $2m'+n'=kn$ and $n' \ge n$, with quotient matrix $kS$.
\end{enumerate}
\end{predln}
\begin{proof}
    Represent a graph $D(km+m',n')$ as a direct product of $m$ copies of graph $D(k,0)$ and graphs $D(a_1,b_1),\ldots,D(a_n,b_n)$, where $2a_i+b_i=k$ for any $i$. Represent the vertices of graph as 
    $(x^1,\ldots,x^m,y^1,\ldots,y^n)$, where $x^i$ is a vertex of $D(k,0)$ and $y^i$ is a vertex of $D(a_i,b_i)$.
    For any $a \in \mathbb{Z}^2_4$ denote $A_a=\{x \in \VV{D(k,0)}:x^{*}_1+\ldots+x^{*}_k=a\}$. 
    Let $B^i_0,\ldots,B^i_3$ be a partition of graph $D(a_i,b_i)$ into disjoint $2$-MDS codes.    

Let $z$ be a vertex of $D(m,n)$. Denote by $C_z=\{(x^1,\ldots,x^m,y^1,\ldots,y^n) \in D(mk+m',n'): x^i \in A_{z^{*}_i}, y^j \in B^j_{z'}\}$.

If vertices $z_1$ and $z_2$ are adjacent in $D(m,n)$ then induced subgraph of $D(m+m',n')$ on the set of vertices $C_{z_1} \cup C_{z_2}$ is $k$-regular bipartite graph, otherwise it is empty graph. Also, $\cup_{z \in \VV{D(m,n)}} C_z=\VV{D(m+m',n')}$ and if $z_1 \ne z_2$  then $C_{z_1}$ and $C_{z_2}$ are disjoint. 
    
Define coloring $f$ of $D(km+m',n')$ as follows: 
if $x \in C_z$ then $f(x)=g(z)$. 

Let $i$ be an arbitrary color and $x$ be an arbitrary vertex from $C_z$ such that $g(z)=i$. The vertex $x$ has neighbours in $C_{z'}$ if and only if $z'$ adjacent to $z$ in $D(m,n)$. Denote set of such neighbours as $N_{z'}$. Let $z_1,\ldots,z_{s_{i,j}}$ are all neighbours of $z$ in $D(m,n)$ with color $j$.  Hence, the set of neighbours of $x$ that has color $j$ equals $N_{z_1} \cup \ldots \cup N_{z_{s_{i,j}}}$.  Since $|N_{z_l}|=k$ for any $l$ we have that $x$ has $ks_{i,j}$ neighbours of color $j$. Therefore, $f$ is a perfect coloring with quotient matrix $kS$.
    
\end{proof}

\begin{lemman}\label{l:d40}
    For arbitrary $a \in \mathbb{Z}^2_4$ the set $C^a=\{x \in \VV{D(4,0)}: x^{*}_1+x^{*}_2+x^{*}_3+x^{*}_4=a\}$ of  vertices of $D(4,0)$ can be partition into codes $C^a_0,\ldots,C^a_3$ with code distance $3$.
\end{lemman}
\begin{proof}
Let $0, 1, \alpha, \alpha^2$ be elements of Galois field $\mathbb{F}_4=GF(4)$.
    Denote a partition of vertices of $D(1,0)$ in the following way: $$A_{00}=\{00,02,20,22\}, A_1=\{01,03,21,23\}, A_{\alpha}=\{10,12,30,32\}, A_{\alpha^2}=\{11,13,31,33\}.$$
    It can be verify directly that if $a \in A_r$ and $b \in A_s$ for any $a,b \in \mathbb{Z}^2_4$ and $r,s \in \mathbb{F}_4$  we have $a+b \in A_{r+s}$. 
    For arbitrary vertex $x$ of $D(4,0)$ denote $y^x=(y^x_1,y^x_2,y^x_3,y^x_4) \in \mathbb{F}^4_4$ where $y^x_i=b$ if $x^{*}_i \in A_b$.
    Define coloring $f:C^{00} \to \{0,1,\alpha,\alpha^2\}$ of vertices from $C^{00}$ in the following way: 
    $f(x)=y^x_2+y^x_3 \alpha + y^x_4 \alpha^2$. 

    Note that $C^{00}$ is an independent set.
    Let $x,z \in C^{00}$ and $d(x,z)=2$. We have that $x$ and $z$ differ in two positions $s,t$. Moreover, $x^{*}_s-z^{*}_s$ and $x^{*}_t-z^{*}_t$ belong to $\{01,03,01,30,11,33\}$ and  $(x^{*}_s-z^{*}_s)+(x^{*}_t-z^{*}_t)=00$. So, the positions of $(y^x-y^z)$ have $2$ zero values and $2$ nonzero equal values. So, $f(x)-f(z)=u(v+w)$, where $u \in \{1,\alpha,\alpha^2\}$ and $v$ and $w$ are different elements from $F_4$. Hence, $f(x)-f(x) \ne 0$. So, any vertices on the distance $2$ between each other have different colors and codes $C^{00}_0=f^{-1}(0)$, $C^{00}_1=f^{-1}(1)$, $C^{00}_2=f^{-1}(\alpha)$ and $C^{00}_3=f^{-1}(\alpha^2)$ form a requirement partition of $C^{00}$. The codes $C^a_0,\ldots,C^a_3$, where $C^a_i=\{x+(a,00,00,00):x \in C^{00}_i\}$, form a requirement partition of $C^a$. 
\end{proof}

\begin{predln}\label{p:splitting}
    Let $g$ be a perfect $k$-coloring in $D(m,n)$ with quotient matrix $S$. Then there is a perfect $4k$-coloring $f$ in $D(4m+c,4n-2c)$, where $c \in\{0,\ldots,2n\}$ with quotient matrix
    \[\displaystyle{
    \begin{pmatrix}
S_{11} & S_{12} & \ldots & S_{1k} \\   
S_{21} & S_{22} & \ldots & S_{2k} \\
\ldots & \ldots & \ldots & \ldots \\
S_{k1} & S_{k2} & \ldots & S_{kk} 
    \end{pmatrix}
    }\]
    where $S_{i,j}$ is $s_{i,j}J_{4 \times 4}$.
\end{predln}
\begin{proof}
As in the Proposition~\ref{p:multiplyk} represent arbitrary vertex $x$ of $D(4m+c,4n-2c)$ as $(x^1,\ldots,x^m,x^{m+1},\ldots,x^{m+n})$, where $x^i$ the vertex of $D(4,0)$ if $i \le m$ and $x^{m+i}$ is a vertex of $D(a_i,b_i)$, where $2a_i+b_i=4$.  For any $a \in \mathbb{Z}^2_4$ denote $C_a=\{x \in \VV{D(4,0)}:x^{*}_1+\ldots+x^{*}_4=a\}$. By Lemma~\ref{l:d40} $C_a$ can be partitioned into codes $C_{(a,0)},\ldots,C_{(a,3)}$ with code distances $3$.
    %Let $C^i_0,\ldots,C^i_3$ be a partition of graph $D(a_i,b_i)$ into $2$-MDS codes. 
    By Proposition~\ref{p:multdiam4} there is a $2$-multipartite perfect coloring in $D(a_i,b_i)$ in colors $(0,0),\ldots,(3,3)$. Denote corresponding color sets as $C^i_{(0,0)},\ldots,C^i_{(3,3)}$.
    
   Define the functions $h:\VV{D(4m+c,4n-2c)} \to \VV{D(m,n)}$ and $r:\VV{D(4m+c,4n-2c)} \to \VV{H(m+n,4)}$ in the following way: 
   $h(x^1,\ldots,x^m,x^{m+1},\ldots,x^{m+n})=z$ and $r(x^1,\ldots,x^m,x^{m+1},\ldots,x^{m+n})=t=(t_1,\ldots,t_{m+n})$, where $x^i \in C_{(z^{*}_i,t_i)}$ if $i \le m$ and $x^{m+j} \in C^{j}_{(z'_j,t_{m+j})}$ otherwise. 

   Let $Y_0,\ldots,Y_3$ be a partition of $H(m+n,4)$ into disjoint $2$-MDS codes. Define $y:\VV{H(m+n,4)} \to \{0,1,2,3\}$ as follows: if $x \in Y_i$ then $s(x)=i$. 
    %Also define function $r:\VV{D(4m+c,4n-2c)} \to \VV{H(m+n,4)}$ in the following way $r(x)=(i_1,\ldots,i_m,j_1,\ldots,j_n)$, where $i_s$ and $j_t$ such that $x^s \in C^s_{i_s}$ and $y^t \in C^t_{}$
    For any $x \in \VV{D(4m+c,4n-c)}$ define coloring $f$ of vertices of $C_z$ as follows $f(x)=(g(h(x)),s(r(x)))$.

Let $x$ be arbitrary vertex in $D(4m+c,4n-c)$ that has color $(i,j)$. 
Denote $h(x)=z$ and $r(x)=(t_1,\ldots,t_{m+n})$.
For vertex $v$ of $D(m,n)$ denote $C_v$ the set of vertices $h^{-1}(v)=\{x \in \VV(D(4m+c,4n-c)):h(x)=v\}$.
Denote $N_v(x)$ the set of neighbours of $x$ in $C_v$.
Note that if $v$ and $z$ are adjacent $N_v(x)$ consists of $4$ vertices with different colors and $N_v(x)$ is empty otherwise. 
Indeed, let $v$ and $z$ differs in one position. Denote this position by $l$.

Suppose, $l>m$. Then $y$ from $N_v(x)$ adjacent to $x$ iff $y^u=x^u$ for any $l \ne u$ and $x^l$ adjacent to $y^l \in C^{l-m}_{(v'_{l-m},b)}$ in $D(a_{l-m},b_{l-m})$ for some $b$. By definition of $2$-multipart perfect coloring for any $b \in \{0,1,2,3\}$ there is exactly one vertex from $C^{l-m}_{(v'_{l-m},b)}$ that adjacent to $x^l$. So, there are $4$ vertices $y_1,\ldots,y_3$ from $N_v$ that adjacent to $x$ and for any $p \ne q$ vertices $r(y_p)$ and $r(y_q)$ from $H(m+n,4)$ differ only in position $l$. By definition of function $s$ for any $p \ne q$ the values $s(r(y_p))$ and $(s(r(y_q)))$ are different. Since $x$ has $4$ neighbours in $N_v$ and they have colors $g(h(x),0),\ldots,(g(h(x)),3)$ once at a time. 
%Hence, for any $i' \in \{1,\ldots,k\}$ the vertex $x$ has $s_{i,i'}$ neighbours of color $(i',i'')$ for any $i''$.

Case $l \le m$ is analogous. A vertex $y$ from $N_v(x)$ adjacent to $x$ iff $y^u=x^u$ for any $l \ne u$ and $x^l$ adjacent to $y^l \in C^{l}_{(v^{*}_{l},b)}$ in $D(4,0)$ for some $b$. 
Since for any $b$ code $C^{l}_{(v^{*}_{l},b)}$ has code distance $3$
for any $b$ there is exactly one vertex from $C^{l}_{(v^{*}_{l},b)}$ that adjacent to $x^l$. So, there are $4$ vertices $y_1,\ldots,y_3$ from $N_v$ that adjacent to $x$ and for any $p \ne q$ vertices $r(y_p)$ and $r(y_q)$ from $H(m+n,4)$ differ only in position $l$. By definition of function $s$ for any $p \ne q$ the values $s(r(y_p))$ and $(s(r(y_q)))$ are different. Since $x$ has $4$ neighbours in $N_v$ and they have colors $g(h(x),0),\ldots,(g(h(x)),3)$ once at a time.  

Therefore, the number of neighbours of $x$ that have color $(i',j')$ equals the number of vertices of $D(m,n)$ that adjacent to $h(x)$ and has a color $i'$. This number equals $s_{i,i'}$. This complete the proof.
\end{proof}

\begin{coroll}
Let there is a perfect $(b,c)$-coloring $g$ in $D(m,n)$.
%$S=\displaystyle{
%    \begin{pmatrix}
%a & b \\   
%c & d
%    \end{pmatrix}}$.
If $b+c>6m+3n$ then there is a perfect $(t_1c+t_2b,(4-t_1)c+(4-t_2)b)$-coloring $f$ in $D(m',n')$, where $2m'+n'=2m+n+(b+c)$ and $t_1,t_2 \in \{0,1,\ldots,4\}$, $t_1+t_2 \ne 0,8$.
\end{coroll}
\begin{proof}
Denote $a=6m+3n-b$ and $d=6m+3n-c$, i.e. $(1,1)$-th and $(2,2)$-th elements of quotient matrix of coloring $g$. 
By Proposition~\ref{p:splitting} there is a perfect coloring in $D(m_1,n_1)$, where $2m_1+n_1=4(2m+n)$, with quotient matrix
\[\displaystyle{
    \begin{pmatrix}
a & \ldots & a & b & \ldots & b \\   
\ldots & \ldots & \ldots & \ldots \\
a & \ldots & a & b & \ldots & b \\
c & \ldots & c & d & \ldots & d \\
\ldots & \ldots & \ldots & \ldots \\
c & \ldots & c & d & \ldots & d
    \end{pmatrix}
    }\].
    By Proposition~\ref{p:latin} there is a perfect $4$-coloring with quotient matrix $(c-a)(J-E)$.
By Theorem~\ref{t:diag} there is a perfect coloring with quotient matrix 
\[\displaystyle{
    \begin{pmatrix}
a & c & c & c & b & b & b & b \\   
c & c & c & c & b & b & b & b \\
c & c & a & c & b & b & b & b \\
c & c & c & a & b & b & b & b \\
c & c & c & c & d & b & b & b \\
c & c & c & c & b & d & b & b \\
c & c & c & c & b & b & d & b \\
c & c & c & c & b & b & b & d
    \end{pmatrix}
    }\].
If we unite first $t_1$ colors and last $t_2$ colors we obtain the coloring $f$ that is perfect with requirement parameters.
\end{proof}

\def\1{node [aa] {\small 1} ++(1,0)}
\def\2{node [bb] {\small 2} ++(1,0)}
\def\3{node [cc] {\small 3} ++(1,0)}
\def\4{node [dd] {\small 4} ++(1,0)}
\def\5{node [ee] {\small 5} ++(1,0)}
\def\6{node [ff] {\small 6} ++(1,0)}
\def\7{node [gg] {\small 7} ++(1,0)}
\def\8{node [hh] {\small 8} ++(1,0)}

\newcommand{\ClipGrid}{\clip [xslant=-0.577] (-1.4,-1.20) rectangle (2.4,2.1);
\draw[xslant=0.577,ystep=.866,xstep=1,draw=black] (-4.9,-2.1) grid (5.4,3.9);
\draw[xslant=-0.577,ystep=9.866,xstep=1,draw=black] (-3.4,-2.1) grid (6.4,3.9);}
\newcommand{\ClipGridV}{\clip  (-1.9,-1.20) rectangle (2.4,2.1);
\draw[xslant=0.577,ystep=.866,xstep=1,draw=black] (-2.6,-2.1) grid (2.6,3.9);
\draw[xslant=-0.577,ystep=9.866,xstep=1,draw=black] (-2.4,-2.1) grid (5.4,3.9);}

\begin{figure}
\centering
\begin{tikzpicture}[scale=0.49,
aa/.style={circle,fill=white,draw=black,inner sep= .2pt},
bb/.style={circle,fill=white,draw=black,inner sep= .2pt},
cc/.style={circle,fill=white,draw=black,inner sep= .2pt},
dd/.style={circle,fill=white,draw=black,inner sep= .2pt},
ee/.style={circle,fill=white,draw=black,inner sep= .2pt},
ff/.style={circle,fill=white,draw=black,inner sep= .2pt},
gg/.style={circle,fill=white,draw=black,inner sep= .2pt},
zz/.style={circle,fill=white,draw=white,inner sep= .2pt},
]
\begin{scope}
\ClipGrid
\draw (-120:1)       \3\4\1\2
++(-4,0) ++(120:1)  \3\4\1\2 
++(-4,0) ++(120:1) \3\4\1\2
++(-4,0) ++(120:1)\3\4\1\2;
\end{scope}
\end{tikzpicture}
\hspace{0.3ex}
\begin{tikzpicture}[scale=0.49,
aa/.style={circle,fill=white,draw=black,inner sep= .2pt},
bb/.style={circle,fill=white,draw=black,inner sep= .2pt},
cc/.style={circle,fill=white,draw=black,inner sep= .2pt},
dd/.style={circle,fill=white,draw=black,inner sep= .2pt},
ee/.style={circle,fill=white,draw=black,inner sep= .2pt},
ff/.style={circle,fill=white,draw=black,inner sep= .2pt},
gg/.style={circle,fill=white,draw=black,inner sep= .2pt},
zz/.style={circle,fill=white,draw=black,inner sep= .2pt},
hh/.style={circle,fill=white,draw=black,inner sep= .2pt},
]
\begin{scope}
\ClipGrid
\draw (-120:1)       \7\8\5\6
++(-4,0) ++(120:1)  \7\8\5\6 
++(-4,0) ++(120:1) \7\8\5\6
++(-4,0) ++(120:1)\7\8\5\6;
\end{scope}
\end{tikzpicture}
\hspace{0.3ex}
\begin{tikzpicture}[scale=0.49,
aa/.style={circle,fill=white,draw=black,inner sep= .2pt},
bb/.style={circle,fill=white,draw=black,inner sep= .2pt},
cc/.style={circle,fill=white,draw=black,inner sep= .2pt},
dd/.style={circle,fill=white,draw=black,inner sep= .2pt},
ee/.style={circle,fill=white,draw=black,inner sep= .2pt},
ff/.style={circle,fill=white,draw=black,inner sep= .2pt},
gg/.style={circle,fill=white,draw=black,inner sep= .2pt},
hh/.style={circle,fill=white,draw=black,inner sep= .2pt},
]
\begin{scope}
\ClipGrid
\draw (-120:1)       \4\3\2\1
++(-4,0) ++(120:1)  \4\3\2\1 
++(-4,0) ++(120:1) \4\3\2\1
++(-4,0) ++(120:1)\4\3\2\1;
\end{scope}
\end{tikzpicture}
\hspace{0.3ex}
\begin{tikzpicture}[scale=0.49,
aa/.style={circle,fill=white,draw=black,inner sep= .2pt},
bb/.style={circle,fill=white,draw=black,inner sep= .2pt},
cc/.style={circle,fill=white,draw=black,inner sep= .2pt},
dd/.style={circle,fill=white,draw=black,inner sep= .2pt},
ee/.style={circle,fill=white,draw=black,inner sep= .2pt},
ff/.style={circle,fill=white,draw=black,inner sep= .2pt},
gg/.style={circle,fill=white,draw=black,inner sep= .2pt},
hh/.style={circle,fill=white,draw=black,inner sep= .2pt},
]
\begin{scope}
\ClipGrid
\draw (-120:1)       \8\7\6\5
++(-4,0) ++(120:1)  \8\7\6\5 
++(-4,0) ++(120:1) \8\7\6\5
++(-4,0) ++(120:1)\8\7\6\5;
\end{scope}
\end{tikzpicture}

\begin{tikzpicture}[scale=0.49,
aa/.style={circle,fill=white,draw=black,inner sep= .2pt},
bb/.style={circle,fill=white,draw=black,inner sep= .2pt},
cc/.style={circle,fill=white,draw=black,inner sep= .2pt},
dd/.style={circle,fill=white,draw=black,inner sep= .2pt},
ee/.style={circle,fill=white,draw=black,inner sep= .2pt},
ff/.style={circle,fill=white,draw=black,inner sep= .2pt},
gg/.style={circle,fill=white,draw=black,inner sep= .2pt},
hh/.style={circle,fill=white,draw=black,inner sep= .2pt},
]
\begin{scope}
\ClipGrid
\draw (-120:1)       \1\2\3\4
++(-4,0) ++(120:1)  \1\2\3\4 
++(-4,0) ++(120:1) \1\2\3\4
++(-4,0) ++(120:1)\1\2\3\4;
\end{scope}
\end{tikzpicture}
\hspace{0.3ex}
\begin{tikzpicture}[scale=0.49,
aa/.style={circle,fill=white,draw=black,inner sep= .2pt},
bb/.style={circle,fill=white,draw=black,inner sep= .2pt},
cc/.style={circle,fill=white,draw=black,inner sep= .2pt},
dd/.style={circle,fill=white,draw=black,inner sep= .2pt},
ee/.style={circle,fill=white,draw=black,inner sep= .2pt},
ff/.style={circle,fill=white,draw=black,inner sep= .2pt},
gg/.style={circle,fill=white,draw=black,inner sep= .2pt},
hh/.style={circle,fill=white,draw=black,inner sep= .2pt},
]
\begin{scope}
\ClipGrid
\draw (-120:1)       \5\6\7\8
++(-4,0) ++(120:1)  \5\6\7\8 
++(-4,0) ++(120:1) \5\6\7\8
++(-4,0) ++(120:1)\5\6\7\8;
\end{scope}
\end{tikzpicture}
\hspace{0.3ex}
\begin{tikzpicture}[scale=0.49,
aa/.style={circle,fill=white,draw=black,inner sep= .2pt},
bb/.style={circle,fill=white,draw=black,inner sep= .2pt},
cc/.style={circle,fill=white,draw=black,inner sep= .2pt},
dd/.style={circle,fill=white,draw=black,inner sep= .2pt},
ee/.style={circle,fill=white,draw=black,inner sep= .2pt},
ff/.style={circle,fill=white,draw=black,inner sep= .2pt},
gg/.style={circle,fill=white,draw=black,inner sep= .2pt},
hh/.style={circle,fill=white,draw=black,inner sep= .2pt},
]
\begin{scope}
\ClipGrid
\draw (-120:1)       \2\1\4\3
++(-4,0) ++(120:1)  \2\1\4\3 
++(-4,0) ++(120:1) \2\1\4\3
++(-4,0) ++(120:1)\2\1\4\3;
\end{scope}
\end{tikzpicture}
\hspace{0.3ex}
\begin{tikzpicture}[scale=0.49,
aa/.style={circle,fill=white,draw=black,inner sep= .2pt},
bb/.style={circle,fill=white,draw=black,inner sep= .2pt},
cc/.style={circle,fill=white,draw=black,inner sep= .2pt},
dd/.style={circle,fill=white,draw=black,inner sep= .2pt},
ee/.style={circle,fill=white,draw=black,inner sep= .2pt},
ff/.style={circle,fill=white,draw=black,inner sep= .2pt},
gg/.style={circle,fill=white,draw=black,inner sep= .2pt},
hh/.style={circle,fill=white,draw=black,inner sep= .2pt},
]
\begin{scope}
\ClipGrid
\draw (-120:1)       \6\5\8\7
++(-4,0) ++(120:1)  \6\5\8\7 
++(-4,0) ++(120:1) \6\5\8\7
++(-4,0) ++(120:1)\6\5\8\7;
\end{scope}
\end{tikzpicture}

\begin{tikzpicture}[scale=0.49,
aa/.style={circle,fill=white,draw=black,inner sep= .2pt},
bb/.style={circle,fill=white,draw=black,inner sep= .2pt},
cc/.style={circle,fill=white,draw=black,inner sep= .2pt},
dd/.style={circle,fill=white,draw=black,inner sep= .2pt},
ee/.style={circle,fill=white,draw=black,inner sep= .2pt},
ff/.style={circle,fill=white,draw=black,inner sep= .2pt},
gg/.style={circle,fill=white,draw=black,inner sep= .2pt},
zz/.style={circle,fill=white,draw=white,inner sep= .2pt},
]
\begin{scope}
\ClipGrid
\draw (-120:1)       \3\4\1\2
++(-4,0) ++(120:1)  \3\4\1\2 
++(-4,0) ++(120:1) \3\4\1\2
++(-4,0) ++(120:1)\3\4\1\2;
\end{scope}
\end{tikzpicture}
\hspace{0.3ex}
\begin{tikzpicture}[scale=0.49,
aa/.style={circle,fill=white,draw=black,inner sep= .2pt},
bb/.style={circle,fill=white,draw=black,inner sep= .2pt},
cc/.style={circle,fill=white,draw=black,inner sep= .2pt},
dd/.style={circle,fill=white,draw=black,inner sep= .2pt},
ee/.style={circle,fill=white,draw=black,inner sep= .2pt},
ff/.style={circle,fill=white,draw=black,inner sep= .2pt},
gg/.style={circle,fill=white,draw=black,inner sep= .2pt},
zz/.style={circle,fill=white,draw=black,inner sep= .2pt},
hh/.style={circle,fill=white,draw=black,inner sep= .2pt},
]
\begin{scope}
\ClipGrid
\draw (-120:1)       \7\8\5\6
++(-4,0) ++(120:1)  \7\8\5\6 
++(-4,0) ++(120:1) \7\8\5\6
++(-4,0) ++(120:1)\7\8\5\6;
\end{scope}
\end{tikzpicture}
\hspace{0.3ex}
\begin{tikzpicture}[scale=0.49,
aa/.style={circle,fill=white,draw=black,inner sep= .2pt},
bb/.style={circle,fill=white,draw=black,inner sep= .2pt},
cc/.style={circle,fill=white,draw=black,inner sep= .2pt},
dd/.style={circle,fill=white,draw=black,inner sep= .2pt},
ee/.style={circle,fill=white,draw=black,inner sep= .2pt},
ff/.style={circle,fill=white,draw=black,inner sep= .2pt},
gg/.style={circle,fill=white,draw=black,inner sep= .2pt},
hh/.style={circle,fill=white,draw=black,inner sep= .2pt},
]
\begin{scope}
\ClipGrid
\draw (-120:1)       \4\3\2\1
++(-4,0) ++(120:1)  \4\3\2\1 
++(-4,0) ++(120:1) \4\3\2\1
++(-4,0) ++(120:1)\4\3\2\1;
\end{scope}
\end{tikzpicture}
\hspace{0.3ex}
\begin{tikzpicture}[scale=0.49,
aa/.style={circle,fill=white,draw=black,inner sep= .2pt},
bb/.style={circle,fill=white,draw=black,inner sep= .2pt},
cc/.style={circle,fill=white,draw=black,inner sep= .2pt},
dd/.style={circle,fill=white,draw=black,inner sep= .2pt},
ee/.style={circle,fill=white,draw=black,inner sep= .2pt},
ff/.style={circle,fill=white,draw=black,inner sep= .2pt},
gg/.style={circle,fill=white,draw=black,inner sep= .2pt},
hh/.style={circle,fill=white,draw=black,inner sep= .2pt},
]
\begin{scope}
\ClipGrid
\draw (-120:1)       \8\7\6\5
++(-4,0) ++(120:1)  \8\7\6\5 
++(-4,0) ++(120:1) \8\7\6\5
++(-4,0) ++(120:1)\8\7\6\5;
\end{scope}
\end{tikzpicture}

\begin{tikzpicture}[scale=0.49,
aa/.style={circle,fill=white,draw=black,inner sep= .2pt},
bb/.style={circle,fill=white,draw=black,inner sep= .2pt},
cc/.style={circle,fill=white,draw=black,inner sep= .2pt},
dd/.style={circle,fill=white,draw=black,inner sep= .2pt},
ee/.style={circle,fill=white,draw=black,inner sep= .2pt},
ff/.style={circle,fill=white,draw=black,inner sep= .2pt},
gg/.style={circle,fill=white,draw=black,inner sep= .2pt},
hh/.style={circle,fill=white,draw=black,inner sep= .2pt},
]
\begin{scope}
\ClipGrid
\draw (-120:1)       \1\2\3\4
++(-4,0) ++(120:1)  \1\2\3\4 
++(-4,0) ++(120:1) \1\2\3\4
++(-4,0) ++(120:1)\1\2\3\4;
\end{scope}
\end{tikzpicture}
\hspace{0.3ex}
\begin{tikzpicture}[scale=0.49,
aa/.style={circle,fill=white,draw=black,inner sep= .2pt},
bb/.style={circle,fill=white,draw=black,inner sep= .2pt},
cc/.style={circle,fill=white,draw=black,inner sep= .2pt},
dd/.style={circle,fill=white,draw=black,inner sep= .2pt},
ee/.style={circle,fill=white,draw=black,inner sep= .2pt},
ff/.style={circle,fill=white,draw=black,inner sep= .2pt},
gg/.style={circle,fill=white,draw=black,inner sep= .2pt},
hh/.style={circle,fill=white,draw=black,inner sep= .2pt},
]
\begin{scope}
\ClipGrid
\draw (-120:1)       \5\6\7\8
++(-4,0) ++(120:1)  \5\6\7\8 
++(-4,0) ++(120:1) \5\6\7\8
++(-4,0) ++(120:1)\5\6\7\8;
\end{scope}
\end{tikzpicture}
\hspace{0.3ex}
\begin{tikzpicture}[scale=0.49,
aa/.style={circle,fill=white,draw=black,inner sep= .2pt},
bb/.style={circle,fill=white,draw=black,inner sep= .2pt},
cc/.style={circle,fill=white,draw=black,inner sep= .2pt},
dd/.style={circle,fill=white,draw=black,inner sep= .2pt},
ee/.style={circle,fill=white,draw=black,inner sep= .2pt},
ff/.style={circle,fill=white,draw=black,inner sep= .2pt},
gg/.style={circle,fill=white,draw=black,inner sep= .2pt},
hh/.style={circle,fill=white,draw=black,inner sep= .2pt},
]
\begin{scope}
\ClipGrid
\draw (-120:1)       \2\1\4\3
++(-4,0) ++(120:1)  \2\1\4\3 
++(-4,0) ++(120:1) \2\1\4\3
++(-4,0) ++(120:1)\2\1\4\3;
\end{scope}
\end{tikzpicture}
\hspace{0.3ex}
\begin{tikzpicture}[scale=0.49,
aa/.style={circle,fill=white,draw=black,inner sep= .2pt},
bb/.style={circle,fill=white,draw=black,inner sep= .2pt},
cc/.style={circle,fill=white,draw=black,inner sep= .2pt},
dd/.style={circle,fill=white,draw=black,inner sep= .2pt},
ee/.style={circle,fill=white,draw=black,inner sep= .2pt},
ff/.style={circle,fill=white,draw=black,inner sep= .2pt},
gg/.style={circle,fill=white,draw=black,inner sep= .2pt},
hh/.style={circle,fill=white,draw=black,inner sep= .2pt},
]
\begin{scope}
\ClipGrid
\draw (-120:1)       \6\5\8\7
++(-4,0) ++(120:1)  \6\5\8\7 
++(-4,0) ++(120:1) \6\5\8\7
++(-4,0) ++(120:1)\6\5\8\7;
\end{scope}
\end{tikzpicture}
\caption{Perfect $8$-coloring of $D(2,0)$ from Lemma~\ref{l:3jd81}.}
\label{f:3jd811}
\end{figure}
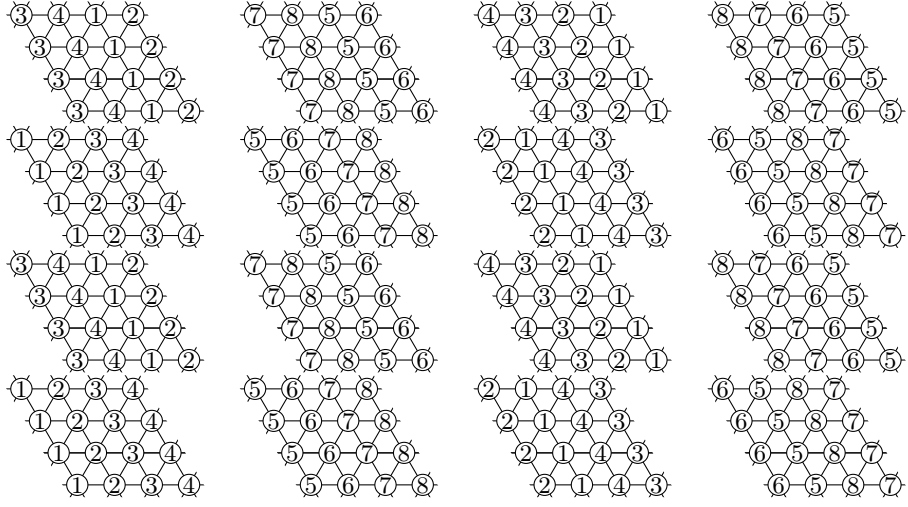

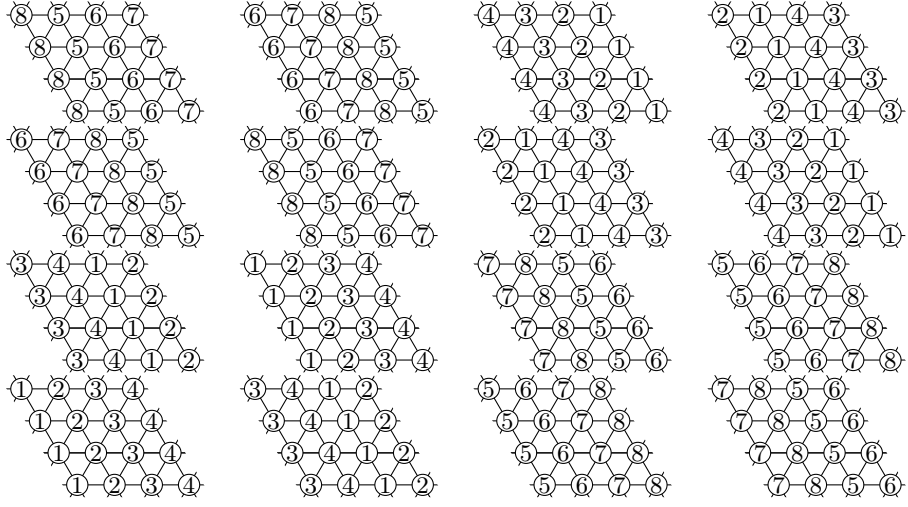
\begin{figure}
\centering
\begin{tikzpicture}[scale=0.49,
aa/.style={circle,fill=white,draw=black,inner sep= .2pt},
bb/.style={circle,fill=white,draw=black,inner sep= .2pt},
cc/.style={circle,fill=white,draw=black,inner sep= .2pt},
dd/.style={circle,fill=white,draw=black,inner sep= .2pt},
ee/.style={circle,fill=white,draw=black,inner sep= .2pt},
ff/.style={circle,fill=white,draw=black,inner sep= .2pt},
gg/.style={circle,fill=white,draw=black,inner sep= .2pt},
hh/.style={circle,fill=white,draw=black,inner sep= .2pt},
]
\begin{scope}
\ClipGrid
\draw (-120:1)       \8\5\6\7
++(-4,0) ++(120:1)  \8\5\6\7 
++(-4,0) ++(120:1) \8\5\6\7
++(-4,0) ++(120:1)\8\5\6\7;
\end{scope}
\end{tikzpicture}
\hspace{0.3ex}
\begin{tikzpicture}[scale=0.49,
aa/.style={circle,fill=white,draw=black,inner sep= .2pt},
bb/.style={circle,fill=white,draw=black,inner sep= .2pt},
cc/.style={circle,fill=white,draw=black,inner sep= .2pt},
dd/.style={circle,fill=white,draw=black,inner sep= .2pt},
ee/.style={circle,fill=white,draw=black,inner sep= .2pt},
ff/.style={circle,fill=white,draw=black,inner sep= .2pt},
gg/.style={circle,fill=white,draw=black,inner sep= .2pt},
zz/.style={circle,fill=white,draw=black,inner sep= .2pt},
hh/.style={circle,fill=white,draw=black,inner sep= .2pt},
]
\begin{scope}
\ClipGrid
\draw (-120:1)       \6\7\8\5
++(-4,0) ++(120:1)  \6\7\8\5 
++(-4,0) ++(120:1) \6\7\8\5
++(-4,0) ++(120:1)\6\7\8\5;
\end{scope}
\end{tikzpicture}
\hspace{0.3ex}
\begin{tikzpicture}[scale=0.49,
aa/.style={circle,fill=white,draw=black,inner sep= .2pt},
bb/.style={circle,fill=white,draw=black,inner sep= .2pt},
cc/.style={circle,fill=white,draw=black,inner sep= .2pt},
dd/.style={circle,fill=white,draw=black,inner sep= .2pt},
ee/.style={circle,fill=white,draw=black,inner sep= .2pt},
ff/.style={circle,fill=white,draw=black,inner sep= .2pt},
gg/.style={circle,fill=white,draw=black,inner sep= .2pt},
hh/.style={circle,fill=white,draw=black,inner sep= .2pt},
]
\begin{scope}
\ClipGrid
\draw (-120:1)       \4\3\2\1
++(-4,0) ++(120:1)  \4\3\2\1 
++(-4,0) ++(120:1) \4\3\2\1
++(-4,0) ++(120:1)\4\3\2\1;
\end{scope}
\end{tikzpicture}
\hspace{0.3ex}
\begin{tikzpicture}[scale=0.49,
aa/.style={circle,fill=white,draw=black,inner sep= .2pt},
bb/.style={circle,fill=white,draw=black,inner sep= .2pt},
cc/.style={circle,fill=white,draw=black,inner sep= .2pt},
dd/.style={circle,fill=white,draw=black,inner sep= .2pt},
ee/.style={circle,fill=white,draw=black,inner sep= .2pt},
ff/.style={circle,fill=white,draw=black,inner sep= .2pt},
gg/.style={circle,fill=white,draw=black,inner sep= .2pt},
hh/.style={circle,fill=white,draw=black,inner sep= .2pt},
]
\begin{scope}
\ClipGrid
\draw (-120:1)       \2\1\4\3
++(-4,0) ++(120:1)  \2\1\4\3 
++(-4,0) ++(120:1) \2\1\4\3
++(-4,0) ++(120:1)\2\1\4\3;
\end{scope}
\end{tikzpicture}

\begin{tikzpicture}[scale=0.49,
aa/.style={circle,fill=white,draw=black,inner sep= .2pt},
bb/.style={circle,fill=white,draw=black,inner sep= .2pt},
cc/.style={circle,fill=white,draw=black,inner sep= .2pt},
dd/.style={circle,fill=white,draw=black,inner sep= .2pt},
ee/.style={circle,fill=white,draw=black,inner sep= .2pt},
ff/.style={circle,fill=white,draw=black,inner sep= .2pt},
gg/.style={circle,fill=white,draw=black,inner sep= .2pt},
hh/.style={circle,fill=white,draw=black,inner sep= .2pt},
]
\begin{scope}
\ClipGrid
\draw (-120:1)       \6\7\8\5
++(-4,0) ++(120:1)  \6\7\8\5 
++(-4,0) ++(120:1) \6\7\8\5
++(-4,0) ++(120:1)\6\7\8\5;
\end{scope}
\end{tikzpicture}
\hspace{0.3ex}
\begin{tikzpicture}[scale=0.49,
aa/.style={circle,fill=white,draw=black,inner sep= .2pt},
bb/.style={circle,fill=white,draw=black,inner sep= .2pt},
cc/.style={circle,fill=white,draw=black,inner sep= .2pt},
dd/.style={circle,fill=white,draw=black,inner sep= .2pt},
ee/.style={circle,fill=white,draw=black,inner sep= .2pt},
ff/.style={circle,fill=white,draw=black,inner sep= .2pt},
gg/.style={circle,fill=white,draw=black,inner sep= .2pt},
hh/.style={circle,fill=white,draw=black,inner sep= .2pt},
]
\begin{scope}
\ClipGrid
\draw (-120:1)       \8\5\6\7
++(-4,0) ++(120:1)  \8\5\6\7 
++(-4,0) ++(120:1) \8\5\6\7
++(-4,0) ++(120:1)\8\5\6\7;
\end{scope}
\end{tikzpicture}
\hspace{0.3ex}
\begin{tikzpicture}[scale=0.49,
aa/.style={circle,fill=white,draw=black,inner sep= .2pt},
bb/.style={circle,fill=white,draw=black,inner sep= .2pt},
cc/.style={circle,fill=white,draw=black,inner sep= .2pt},
dd/.style={circle,fill=white,draw=black,inner sep= .2pt},
ee/.style={circle,fill=white,draw=black,inner sep= .2pt},
ff/.style={circle,fill=white,draw=black,inner sep= .2pt},
gg/.style={circle,fill=white,draw=black,inner sep= .2pt},
hh/.style={circle,fill=white,draw=black,inner sep= .2pt},
]
\begin{scope}
\ClipGrid
\draw (-120:1)       \2\1\4\3
++(-4,0) ++(120:1)  \2\1\4\3 
++(-4,0) ++(120:1) \2\1\4\3
++(-4,0) ++(120:1)\2\1\4\3;
\end{scope}
\end{tikzpicture}
\hspace{0.3ex}
\begin{tikzpicture}[scale=0.49,
aa/.style={circle,fill=white,draw=black,inner sep= .2pt},
bb/.style={circle,fill=white,draw=black,inner sep= .2pt},
cc/.style={circle,fill=white,draw=black,inner sep= .2pt},
dd/.style={circle,fill=white,draw=black,inner sep= .2pt},
ee/.style={circle,fill=white,draw=black,inner sep= .2pt},
ff/.style={circle,fill=white,draw=black,inner sep= .2pt},
gg/.style={circle,fill=white,draw=black,inner sep= .2pt},
hh/.style={circle,fill=white,draw=black,inner sep= .2pt},
]
\begin{scope}
\ClipGrid
\draw (-120:1)       \4\3\2\1
++(-4,0) ++(120:1)  \4\3\2\1 
++(-4,0) ++(120:1) \4\3\2\1
++(-4,0) ++(120:1)\4\3\2\1;
\end{scope}
\end{tikzpicture}

\begin{tikzpicture}[scale=0.49,
aa/.style={circle,fill=white,draw=black,inner sep= .2pt},
bb/.style={circle,fill=white,draw=black,inner sep= .2pt},
cc/.style={circle,fill=white,draw=black,inner sep= .2pt},
dd/.style={circle,fill=white,draw=black,inner sep= .2pt},
ee/.style={circle,fill=white,draw=black,inner sep= .2pt},
ff/.style={circle,fill=white,draw=black,inner sep= .2pt},
gg/.style={circle,fill=white,draw=black,inner sep= .2pt},
zz/.style={circle,fill=white,draw=white,inner sep= .2pt},
]
\begin{scope}
\ClipGrid
\draw (-120:1)       \3\4\1\2
++(-4,0) ++(120:1)  \3\4\1\2 
++(-4,0) ++(120:1) \3\4\1\2
++(-4,0) ++(120:1)\3\4\1\2;
\end{scope}
\end{tikzpicture}
\hspace{0.3ex}
\begin{tikzpicture}[scale=0.49,
aa/.style={circle,fill=white,draw=black,inner sep= .2pt},
bb/.style={circle,fill=white,draw=black,inner sep= .2pt},
cc/.style={circle,fill=white,draw=black,inner sep= .2pt},
dd/.style={circle,fill=white,draw=black,inner sep= .2pt},
ee/.style={circle,fill=white,draw=black,inner sep= .2pt},
ff/.style={circle,fill=white,draw=black,inner sep= .2pt},
gg/.style={circle,fill=white,draw=black,inner sep= .2pt},
zz/.style={circle,fill=white,draw=black,inner sep= .2pt},
hh/.style={circle,fill=white,draw=black,inner sep= .2pt},
]
\begin{scope}
\ClipGrid
\draw (-120:1)       \1\2\3\4
++(-4,0) ++(120:1)  \1\2\3\4 
++(-4,0) ++(120:1) \1\2\3\4
++(-4,0) ++(120:1)\1\2\3\4;
\end{scope}
\end{tikzpicture}
\hspace{0.3ex}
\begin{tikzpicture}[scale=0.49,
aa/.style={circle,fill=white,draw=black,inner sep= .2pt},
bb/.style={circle,fill=white,draw=black,inner sep= .2pt},
cc/.style={circle,fill=white,draw=black,inner sep= .2pt},
dd/.style={circle,fill=white,draw=black,inner sep= .2pt},
ee/.style={circle,fill=white,draw=black,inner sep= .2pt},
ff/.style={circle,fill=white,draw=black,inner sep= .2pt},
gg/.style={circle,fill=white,draw=black,inner sep= .2pt},
hh/.style={circle,fill=white,draw=black,inner sep= .2pt},
]
\begin{scope}
\ClipGrid
\draw (-120:1)       \7\8\5\6
++(-4,0) ++(120:1)  \7\8\5\6 
++(-4,0) ++(120:1) \7\8\5\6
++(-4,0) ++(120:1)\7\8\5\6;
\end{scope}
\end{tikzpicture}
\hspace{0.3ex}
\begin{tikzpicture}[scale=0.49,
aa/.style={circle,fill=white,draw=black,inner sep= .2pt},
bb/.style={circle,fill=white,draw=black,inner sep= .2pt},
cc/.style={circle,fill=white,draw=black,inner sep= .2pt},
dd/.style={circle,fill=white,draw=black,inner sep= .2pt},
ee/.style={circle,fill=white,draw=black,inner sep= .2pt},
ff/.style={circle,fill=white,draw=black,inner sep= .2pt},
gg/.style={circle,fill=white,draw=black,inner sep= .2pt},
hh/.style={circle,fill=white,draw=black,inner sep= .2pt},
]
\begin{scope}
\ClipGrid
\draw (-120:1)       \5\6\7\8
++(-4,0) ++(120:1)  \5\6\7\8 
++(-4,0) ++(120:1) \5\6\7\8
++(-4,0) ++(120:1)\5\6\7\8;
\end{scope}
\end{tikzpicture}

\begin{tikzpicture}[scale=0.49,
aa/.style={circle,fill=white,draw=black,inner sep= .2pt},
bb/.style={circle,fill=white,draw=black,inner sep= .2pt},
cc/.style={circle,fill=white,draw=black,inner sep= .2pt},
dd/.style={circle,fill=white,draw=black,inner sep= .2pt},
ee/.style={circle,fill=white,draw=black,inner sep= .2pt},
ff/.style={circle,fill=white,draw=black,inner sep= .2pt},
gg/.style={circle,fill=white,draw=black,inner sep= .2pt},
hh/.style={circle,fill=white,draw=black,inner sep= .2pt},
]
\begin{scope}
\ClipGrid
\draw (-120:1)       \1\2\3\4
++(-4,0) ++(120:1)  \1\2\3\4 
++(-4,0) ++(120:1) \1\2\3\4
++(-4,0) ++(120:1)\1\2\3\4;
\end{scope}
\end{tikzpicture}
\hspace{0.3ex}
\begin{tikzpicture}[scale=0.49,
aa/.style={circle,fill=white,draw=black,inner sep= .2pt},
bb/.style={circle,fill=white,draw=black,inner sep= .2pt},
cc/.style={circle,fill=white,draw=black,inner sep= .2pt},
dd/.style={circle,fill=white,draw=black,inner sep= .2pt},
ee/.style={circle,fill=white,draw=black,inner sep= .2pt},
ff/.style={circle,fill=white,draw=black,inner sep= .2pt},
gg/.style={circle,fill=white,draw=black,inner sep= .2pt},
hh/.style={circle,fill=white,draw=black,inner sep= .2pt},
]
\begin{scope}
\ClipGrid
\draw (-120:1)       \3\4\1\2
++(-4,0) ++(120:1)  \3\4\1\2 
++(-4,0) ++(120:1) \3\4\1\2
++(-4,0) ++(120:1)\3\4\1\2;
\end{scope}
\end{tikzpicture}
\hspace{0.3ex}
\begin{tikzpicture}[scale=0.49,
aa/.style={circle,fill=white,draw=black,inner sep= .2pt},
bb/.style={circle,fill=white,draw=black,inner sep= .2pt},
cc/.style={circle,fill=white,draw=black,inner sep= .2pt},
dd/.style={circle,fill=white,draw=black,inner sep= .2pt},
ee/.style={circle,fill=white,draw=black,inner sep= .2pt},
ff/.style={circle,fill=white,draw=black,inner sep= .2pt},
gg/.style={circle,fill=white,draw=black,inner sep= .2pt},
hh/.style={circle,fill=white,draw=black,inner sep= .2pt},
]
\begin{scope}
\ClipGrid
\draw (-120:1)       \5\6\7\8
++(-4,0) ++(120:1)  \5\6\7\8 
++(-4,0) ++(120:1) \5\6\7\8
++(-4,0) ++(120:1)\5\6\7\8;
\end{scope}
\end{tikzpicture}
\hspace{0.3ex}
\begin{tikzpicture}[scale=0.49,
aa/.style={circle,fill=white,draw=black,inner sep= .2pt},
bb/.style={circle,fill=white,draw=black,inner sep= .2pt},
cc/.style={circle,fill=white,draw=black,inner sep= .2pt},
dd/.style={circle,fill=white,draw=black,inner sep= .2pt},
ee/.style={circle,fill=white,draw=black,inner sep= .2pt},
ff/.style={circle,fill=white,draw=black,inner sep= .2pt},
gg/.style={circle,fill=white,draw=black,inner sep= .2pt},
hh/.style={circle,fill=white,draw=black,inner sep= .2pt},
]
\begin{scope}
\ClipGrid
\draw (-120:1)       \7\8\5\6
++(-4,0) ++(120:1)  \7\8\5\6 
++(-4,0) ++(120:1) \7\8\5\6
++(-4,0) ++(120:1)\7\8\5\6;
\end{scope}
\end{tikzpicture}
\caption{Perfect $8$-colorings of $D(1,2)$ from Lemma~\ref{l:3jd82}.}
\label{f:3jd812}
\end{figure}

\begin{lemman}\label{l:3jd81}
There is a perfect coloring in $D(2,0)$ and $D(1,2)$ with quotient matrix
\[\displaystyle{
    \begin{pmatrix}
2 & \ldots & 2 & 1 & \ldots & 1 \\   
\ldots & \ldots & \ldots  & \ldots & \ldots &  \ldots \\
2 & \ldots & 2 & 1 & \ldots & 1 \\ 
1 & \ldots & 1 & 2 & \ldots & 2 \\ 
\ldots & \ldots & \ldots  & \ldots & \ldots &  \ldots \\
1 & \ldots & 1 & 2 & \ldots & 2 
    \end{pmatrix}
    }\].
\end{lemman}
\begin{proof}
These colorings are shown on Figures~\ref{f:3jd811}-~\ref{f:3jd812}. 
\end{proof}

\begin{lemman}\label{l:3jd82}
There is a perfect $8$-coloring in $D(m,n)$, where $2m+n=4$, with quotient matrix
\[\displaystyle{
    \begin{pmatrix}
1 & \ldots & 1 & 2 & \ldots & 2 \\   
\ldots & \ldots & \ldots  & \ldots & \ldots &  \ldots \\
1 & \ldots & 1 & 2 & \ldots & 2 \\ 
2 & \ldots & 2 & 1 & \ldots & 1 \\ 
\ldots & \ldots & \ldots  & \ldots & \ldots &  \ldots \\
2 & \ldots & 2 & 1 & \ldots & 1 
    \end{pmatrix}
    }\].
\end{lemman}
\begin{proof}
By Proposition~\ref{p:multdiam4} there is a $2$-multipartite perfect coloring in colors $(0,0),\ldots,(3,3)$. If we unit colors $(0,j)$ and $(1,j)$ and unit colors $(2,j)$ and $(3,j)$ for any $j$ we obtain requirement coloring.
\end{proof}

\begin{predln}\label{p:3j8}
There is a perfect $8$-coloring in $D(m,n)$, where $2m+n=8$, with quotient matrix $3J$. 
\end{predln}
\begin{proof}
Represent colors as $(i,j)$, where $i \in \{0,1\}$, $j \in \{0,1,2,3\}$.
Represent $D(m,n)$ as a direct product of graphs $D(m_1,n_1)$ and $D(m_2,n_2)$, where $2m_1+n_1=2m_2+n_2=4$, $m_2 \ge 1$.  Represent vertices of $D(m,n)$ as $(x,y)$, where $x \in \VV{D(m_1,n_1)}$, $y \in \VV{D(m_2,n_2)}$.
Let $g_1$ be a perfect coloring of $D(m_1,n_1)$ with quotient matrix from Lemma~\ref{l:3jd82}. Let $g_2$ be a perfect coloring of $D(m_2,n_2)$ with quotient matrix from Lemma~\ref{l:3jd81}. Denote coloring $f$ of $D(m,n)$ in the following way:
$$f(x,y)=g_1(x)+g_2(y).$$

Let $(x,y)$ be a vertex of color $(i,j)$. Denote by $N_1(x,y)$ the number of vertices $(x,y')$ that adjacent to $(x,y)$. Denote by $N_2(x,y)$ the number of vertices $(x',y)$ that adjacent to $(x,y)$. The neighbourhood $N(x,y)$ of $(x,y)$ equals $N_1(x,y) \cup N_2(x,y)$.
Let us count the number of vertices of $(r,s)$ that adjacent to $(x,y)$.
The number of such vertices in $N_1(x,y)$ equals $2$ if $i=r$ and $1$ if $i \ne r$. The number of such vertices in $N_2(x,y)$ equals $1$ if $i=r$ and $2$ if $i \ne r$. So, this numbers equals $3$ in any case. 
\end{proof}

\section{Partitions of 2-MDS codes}

\begin{predln}\label{p:mdspart}
     If $n=2^k$ for some $k \ge 2$ then there is a $2$-MDS code in $H(n,4)$ that can be partitioned into $2^k$ codes with distance at least $3$.
\end{predln}
\begin{proof}
    Consider a Galois field $\mathbb{F}_{2^k}=GF(2^k)$ of order $2^k$, with elements $0,1,\alpha,\ldots,\alpha^{2^k-2}$ for some primitive element $\alpha$.
    The elements of this field can be treated as polynomials of degree less than $k$ with coefficient from $\{0,1\}$. 
   
    Consider graph $H(2^k,2^k)$. Define $C'$ as a linear code with check matrix 
\begin{equation*}
H=\begin{pmatrix}
1 & 1 & 1 & 1 & \ldots & 1 & 1 
%\\0 & 1 & \alpha & \alpha^2 & \ldots & \alpha^{p^k-3} & \alpha^{p^k-2} 
\end{pmatrix}
\end{equation*}
These codes can be partition into $2^k$ codes, $R_1,\ldots,R_{2^k}$, with code distance $3$ as a cosets of linear code with check matrix 
\begin{equation*}
H=\begin{pmatrix}
1 & 1 & 1 & 1 & \ldots & 1 & 1 \\
0 & 1 & \alpha & \alpha^2 & \ldots & \alpha^{2^k-3} & \alpha^{2^k-2} 
\end{pmatrix}
\end{equation*}
where $R_1$ corresponds to syndrome $(0,0)$, $R_2$ corresponds to syndrome $(0,1)$ and $R_i$, $i \in \{3,\ldots,2^k\}$ corresponds to syndrome $(0,\alpha^{i-2})$.

Denote by $Q=\{0,1,x,x+1\}$ the set of elements of $\mathbb{F}_{2^k}$ that are polynomials of degree less that $2$.
Consider the set of vertices 
$$A=\{(x_1,\ldots,x_n) \in \VV{H(2^k,2^k)}: x_i \in Q)\}.$$  

The induced subgraph on the set $A$ isomorphic to the Hamming graph $H(n,4)$.
Denote this subgraph by $A(n,4)$ and denote by $d_A(x,y)$ the distance between vertices $x$ and $y$ in graph $A(n,4)$. 
We have that for any vertices $x$ and $y$ from $A$ the distances $d(x,y)$ and $d_A(x,y)$ are equal.

Denote $C=C' \cap A$. A code $C$ is a $2$-MDS code in $A(n,4)$. Indeed, if we fixed values from $Q$ in any $n-1$ positions then there is exactly one element from $Q$ such that sum of values in all $n$ positions equals $0$.   

Denote $L_i=R_i \cap A$. 
Since $(0,\ldots,0)$ and $(1,0,\ldots,1,\ldots,0,\ldots,0)$, where $1$ on the $i$-th position, belong to $C$, then $L_i$ is not empty for any $i$. 
Also for any $i$ we have that if $x$ and $y$ are vertices from $L_i$ then $d_A(x,y)=d(x,y)$ not less than $3$.
So, codes $L_1,\ldots,L_{2^k}$ is a requirement partition.
\end{proof} 

\begin{remarkn}\label{r:mdspart}
    Proposition~\ref{p:mdspart} also is true for any Hamming graphs $H(p^k,p^m)$, where $k \ge m$ and $p$ is prime. A partition into $p^k$ codes with code distance at least $3$ can be constructed  as cosets of additive completely regular code with intersection array $[p^k(p^m-1),1;p^k-1,p^k(p^m-1)]$. These completely regular codes were mentioned in \cite{BRZ:crg}, but only for case when codes are linear which implies that $k$ is divisible by $m$. In case when $k$ is not divisible by $m$
    such code can be constructed as dual to known additive two-weight codes in $H(n,q)$ with parameters: length $n=p^k$, cardinality $N=p^{k+m}$ and code distance $p^{k-m}(p^m-1)$ (see in \cite{SZZ:equidistance}, \cite{BZZ:cover2}). These two-codes constructed with using of difference matrices. 
    %Completely regular code with intersection array $[p^k(p^m-1),1;p^k-1,p^k(p^m-1)]$ were mentioned in \cite{BRZ:crg}, but only for case $k$ is divisible by $m$.
    %It worth to mention that such completely regular code is Orthogonal array with parameters $(p^{km-k-m},p^k,p^m,p^{k-m}(q^m-1)-1)$ that has minimum possible index (because punctured code lies on Bierbrauer-Friedman bound).
\end{remarkn}

\begin{coroll}\label{c:gammamdspartcol}
Let $2m+n=\gamma \cdot 2^k$ some positive integer $\gamma$ and $k$ such that $k \ge 2$ and if $m>0$ then $k \ne 3$. Then there is a perfect 
$(2^k+1)$-colorings $f$ in $    D(m,n)$ with  quotient matrix $S$ equals
    \[\displaystyle{\begin{pmatrix}
0 & 0 & 0 & \ldots & 0 & \gamma(6m+3n)  \\
0 & 0 & 0 & \ldots & 0 & \gamma(6m+3n) \\
\ldots & \ldots & \ldots & \ldots & \ldots \\
0 & 0 & 0 & \ldots & 0 & \gamma(6m+3n) \\
\gamma & \gamma & \gamma & \ldots & \gamma & \gamma(4m+2n) 
\end{pmatrix}}\]
\end{coroll}
\begin{proof}
Let $\gamma=1$ and $m=0$.
    By Proposition~\ref{p:mdspart} there is a $2$-MDS code $C$ that can be partitioned into codes $L_1,\ldots,L_{2^k}$ with code distance $3$. Let $x$ be arbitrary vertex that does not belong to $C$. This vertex adjacent to $2^k$ vertices from $C$ and it can't have more than one neighbour in $L_i$ for any $i$, otherwise we have a contradiction with code distance. 
    Hence, the coloring , where $f(x)=i$ if $x \in L_i$ for $i=1,\ldots,2^k$ and $f(x)=2^k+1$ if $x$ does not belong to $C$, is required perfect coloring.

In \cite{BKMTV:perfcolinham} it was proven that if there is a perfect coloring in $H(n,q)$ with quotient matrix $S$ then there is a perfect coloring in $H(\gamma n,q)$ with quotient matrix $\gamma S$. 

If $m \ne 0$ we can apply Proposition~\ref{p:splitting} to same coloring in $H(\gamma \cdot 2^{k-2},4)$ and unit last four colors. For $k=2$ the statement follows from Proposition~\ref{p:multdiam4}.
\end{proof}

\begin{coroll}\label{c:multipart}
For any positive integer $k$ such that $k \ge 2$ and if $n \ne 0$ then $k \ne 3$ there is a $k$-multipartite perfect coloring in $D(m,n)$, where $2m+n=2^k$.
\end{coroll}
\begin{proof}
Let $m=0$. 
By Proposition~\ref{p:mdspart} there is a $2$-MDS code $C^0$ that can be partitioned into codes $L^0_1,\ldots,L^0_{2^k}$ with code distance $3$. 
Denote by $e_a$ the vertex $(a,0,\ldots,0)$ where $a \in \{1,2,3\}$. The code  $C^a=\{x+e_a: x \in C^0\}$ is a $2$-MDS code and the codes 
$L^a_1,\ldots,L^a_{2^k}$, where
$L^a_i=\{x+e_a:x \in L^0_i\}$, partition $C^a$ into codes with code distance $3$.
The codes $L^0_1,\ldots,L^3_{2^k}$ are satisfy the condition of Lemma~\ref{l:multipart} which implies the statement.

If $m>0$ we can apply Proposition~\ref{p:splitting} to $(k-1)$-multipartite perfect coloring in $H(2^{k-1},4)$. If $k=2$ the statement follows from Lemma~\ref{l:multipart} and Proposition~\ref{p:multdiam4}.
\end{proof}

\begin{coroll}\label{c:3jham}
We have the following.
\begin{enumerate}
\item Let $2m+n=2^k$ for some positive integer $k \ge 2$. Then there is $2^k$-coloring in $D(m,n)$ with quotient matrix $3J$.
\item Let $2m+n=2^k-1$ for some positive $k$ such that $k \ge 2$ and if $m>0$ then $k \ne 3$. Then there is a perfect $2^k$-coloring in $D(m,n)$ with quotient matrix $3(J-E)$.
\end{enumerate}
\end{coroll}
\begin{proof}
1) This is follows from Corollary~\ref{c:multipart}, Lemma~\ref{l:3jmult} and Proposition~\ref{p:3j8}.

2) Note that in this case $n \ne 0$. By Corollary~\ref{c:gammamdspartcol} there is a $2$-MDS code $L^0$ that can be partitioned into codes $L^0_1,\ldots,L^0_{2^l}$ such that any vertex that not belong to $L^0$ adjacent to exactly one vertex from $L^0_j$ for any $j$.
Denote by $e_a$ the vertex $(00,\ldots,00;0,\ldots,0,a)$ of $D(m,n)$ where $a \in \{1,2,3\}$. The code  $L^a=\{x+e_a: x \in L^0\}$ is a $2$-MDS code and the codes 
$L^0_1,\ldots,L^3_{2^k}$ form partition that correspond to $k$-multipartite perfect coloring. The coloring $f$, where $f(x)=i$ if $x \in L^0_{i} \cup \ldots \cup L^3_i$, it is a perfect $2^k$-coloring with quotient matrix $3J$. Moreover, any two vertices that differ in the last position have the same color. So, coloring $f'$ such that $f'(x^{*}_1,\ldots,x^{*}_m;x'_1,\ldots,x'_{n-1})=f(x^{*}_1,\ldots,x^{*}_m;x'_1,\ldots,x'_{n-1},0)$ of $D(m,n-1)$ is a perfect with quotient matrix $3(J-E)$. 
\end{proof}

%\begin{coroll}\label{c:gammamdspartcoldoob}
%If $2m+n=\gamma \cdot 2^k$ for some positive integer $k \ge 4$ and $\gamma$, %then there is an $(2^k+1)$-perfect colorings $f$ in $D(m,n)$ with  quotient %matrix $S$ equals
%    \[\displaystyle{\begin{pmatrix}
%0 & 0 & 0 & \ldots & 0 & \gamma(6m+3n)  \\
%0 & 0 & 0 & \ldots & 0 & \gamma(6m+3n) \\
%\ldots & \ldots & \ldots & \ldots & \ldots \\
%0 & 0 & 0 & \ldots & 0 & \gamma(6m+3n) \\
%\gamma & \gamma & \gamma & \ldots & \gamma & \gamma(4m+2n) 
%\end{pmatrix}}\]
%\end{coroll}
%\begin{proof}
%If we apply Proposition~\ref{p:splitting} to coloring from %Corollary~\ref{c:gammamdspartcol} and unit the last $4$ colors we obtain the %requirement coloring. 
%\end{proof}

%\begin{coroll}\label{c:multipartite}
%If $2m+n=2^k$, where $k \le 2$ and if $k=3$ then $m=0$, there is $2^{k+2}$ %perfect coloring with colors $(i,j)$, $i\in\{0,\ldots,3\}$, $j \in \
%{0,1,\ldots,2^k-1\}$, such that for any $i,j,r,s$ any vertex of color %$(i,j)$ has $1$ neighbour of color $(r,s)$ if $r \ne i$ and $0$ neighbours %if $r=i$. 
%\end{coroll}
%\begin{proof}
%Ai in the proof of Corollary
%\end{proof}

\begin{coroll}\label{c:rad2}
    Let $2m+n=\gamma \cdot 2^k$ for some odd $\gamma$ and $k \ge 4$. Then there is a completely regular code $C$ in $D(m,n)$ with quotient matrix equals
    
     \[\displaystyle{\begin{pmatrix}
0 \ \ \  & 6m+3n & 0 \\
c \ \ \  & 4m+2n & 2m+n-c \\ 
0 \ \ \ & 6m+3n & 0
\end{pmatrix}}\]
if and only if $2m+n$ divides by $4$, $c=\gamma \cdot b$ for some $b \in \{1,2^k-1\}$.
\end{coroll}
\begin{proof}
Let $f$ be a perfect coloring that corresponds to distance coloring of $C$.
From the parameters we have $(6m+3n)|f^{-1}(0)|=c|f^{-1}(1)|$ and $(6m+3n)|f^{-1}(2)|=(2m+n-c)|f^{-1}(1)|$. Hence, the set $f^{-1}(0) \cup f^{-1}(2)$ is independence and has cardinality $4^{2m+n-1}$, so, it is a $2$-MDS code. 
    Also, we have $|f^{-1}(0)|=\displaystyle{\frac{c \cdot 4^{2m+n-1}}{2m+n}}$.
    %Let $2m+n=\gamma \cdot 2^k$, where $\gamma$ is odd. 
    Since $|f^{-1}(0)|$ is integer then $c=b \cdot \gamma$ for some integer $b \in \{1,\ldots,2^k-1\}$. 
    On the other hand, by Corollary~\ref{c:gammamdspartcol} there is a perfect coloring in $D(m,n)$ with quotient matrix equals \[\displaystyle{\begin{pmatrix}
0 & 0 & 0 & \ldots & 0 & \gamma(6m+3n)  \\
0 & 0 & 0 & \ldots & 0 & \gamma(6m+3n) \\
\ldots & \ldots & \ldots & \ldots & \ldots \\
0 & 0 & 0 & \ldots & 0 & \gamma(6m+3n) \\
\gamma & \gamma & \gamma & \ldots & \gamma & \gamma(4m+2n) 
\end{pmatrix}}\]
The union of the first $b$ colors of this coloring is required completely regular code.
\end{proof}

%\begin{coroll}\label{c:3jcol}
%\begin{enumerate}
%    \item For any $l \ge 2$, $l \ne 3$ there is a $2^l$-perfect coloring in %$D(m,0)$, where $2m=2^l$, with quotient matrix $3J$.
%    \item For any $l \ge 2$. $l \ne 3$ there is a $2^l$-perfect coloring in %$D(m,n)$, where $2m+n=2^l-1$, with quotient matrix $3(J-E)$. 
%\end{enumerate}
%\end{coroll}
%\begin{proof}
%By Corollary~\ref{c:3jham} there is a $2^l$-perfect coloring in $H(2^l-1,4)$ %with quotient matrix $3(J-E)$. Hence there is a $2^l$-perfect coloring in %$H(2^l,4)$. Also there is a perfect coloring in $D(1,1)$ and $D(0,3)$ with
%quotient matrix $3(J_{4 \times 4}-E_{4 \times 4})$ (a partition into %disjoint $2$-MDS codes).
%Therefore, the requirement coloring can be obtained by %Proposition~\ref{p:splitting} and Theorem~\ref{t:diag}.
%\end{proof}

\section{Existence of perfect (b,c)-colorings}

\subsection{Perfect $2^l$-colorings}

\begin{predln}\label{p:2lham}
We have the following.
\begin{enumerate}
    \item For any positive integer $l$ and $s$ there is a perfect 
    $4^l$-coloring in $D(0,n)$, where $n=\displaystyle{\frac{s(4^l-1)}{3}}$, with quotient matrix $s(J-E)$.
    \item For any positive integer $l \ge 2$ and $s \ge 2$ there is a perfect
    $2^{2l-1}$-coloring in $D(0,n)$, where $n=\displaystyle{\frac{s(2^{2l-1}-1)+a}{3}}$, with quotient matrix $s(J-E)+aE$, where $a=0$ if ($s=0 \mod 3$), $a=1$ if ($s=2 \mod 3$) and $a=2$ if ($s=1 \mod 3$). 
    \item For any positive integer $l \ge 2$ there is a perfect $2^{2l-1}$-coloring in $D(m,n)$, where $2m+n=\displaystyle{\frac{(2^{2l}-1)}{3}}$, with quotient matrix $2(J-E)+E$ 
    \item For any positive integer $l \ge 2$ there is a perfect $2^{2l-1}$-coloring in $D(m,0)$, where $2m=\displaystyle{\frac{2(2^{2l-1}+1)}{3}}$, with quotient matrix $2(J-E)+4E$ 
\end{enumerate}
\end{predln}
\begin{proof}
1) By Corollary~\ref{c:JE} there is a perfect $4^l$-coloring in $H(\frac{4^l-1}{3},4)$ with quotient matrix $(J-E)$. By Proposition~\ref{p:multiplyk} there is a perfect $4^l$-coloring in $H(\frac{s(4^l-1)}{3},4)$ with quotient matrix $s(J-E)$.

2) By Corollary~\ref{c:3jham} there is a $2^{2l-1}$ perfect coloring in $H(2^{2l-1}-1,4)$ with quotient matrix $3(J-E)$. Also there is a partition into $2^{2l-1}$ disjoint $2$-fold perfect codes in $H(\frac{4^l-1}{3},4)$. This corresponds to perfect $2^{2l-1}$-coloring with quotient matrix $2(J-E)+E$. If $s=3k+1$ by Theorem~\ref{t:diag} there is a perfect $2^{2l-1}$-coloring with quotient matrix 
$3(J-E)+\ldots+3(J-E)+(2(J-E)+E)+(2(J-E)+E)$ that equals $(3k+1)(J-E)+2E$.
If $s=3k$ and $s=3k+1$ construction is analogous.

3) This coloring can be obtained using the partition into $2$-fold perfect codes in $D(m,n)$, where $2m+n=\frac{2^{2l}-1}{3}$.

4) Let us prove the statement by induction on $l$.
Let $l=1$. There is a perfect $(2,2)$-coloring (the set $\{00,01,02,03,10,11,12,13\}$ forms the first color) in $D(1,0)$.
Suppose there is a perfect $2^{2l-1}$-coloring $f$ with requirement coloring.
By Corollary~\ref{c:multipart} there is a $2l$-multipartite perfect coloring  in $D(m,n)$, where $2m=2^{2l}$. If we unit colors $(i,2j-1)$ and $(i,2j)$ for any $i$ and $j$ and apply Theorem~\ref{t:diag} to obtained perfect coloring and coloring $f$ (taken $4$ times) we obtain the requirement perfect $2^{2l+1}$-coloring.
\end{proof}

%\begin{theoreman}
 %   \begin{enumerate}
 %       \item For any $l$ and even $s$ there is $4^l$-perfect coloring in %$D(m,n)$, where $2m+n=\frac{s(4^l-1)}{3}$, with quotient matrix $s(J-E)$. 
 %       \item For any $l$ and odd $s$ there is $4^l$-perfect coloring in %$D(m,n)$, where $2m+n=\frac{s(4^l-1)}{3}$, with quotient matrix $s(J-E)$.
 %       \item For any $l$ and odd $s$ there is $4^l$-perfect coloring in %$D(m,0)$, where $2m=\frac{s(4^l-1)+3}{3}$, with quotient matrix $s(J-E)+3E$.
%\end{enumerate}
%\end{theoreman}

\begin{predln}\label{p:4lcol}
We have the following.
    \begin{enumerate}
        \item For any $l$ and $s$ there is perfect $4^l$-coloring in $D(m,n)$, where $\displaystyle{2m+n=\frac{s(4^l-1)}{3}}$, with quotient matrix $s(J-E)$. 
        %\item For any $l$ and odd $s$ there is perfect $4^l$-coloring in $D(m,n)$, where $2m+n=\frac{s(4^l-1)}{3}$, with quotient matrix $s(J-E)$.
        \item For any $l$ and odd $s \ge 3$ there is perfect $4^l$-coloring in $D(m,0)$, where $\displaystyle{2m=\frac{s(4^l-1)+3}{3}}$, with quotient matrix $s(J-E)+3E$.
\end{enumerate}
\end{predln}
\begin{proof}
    By Proposition~\ref{p:2lham} there is a perfect $4^l$-coloring in $D(0,n)$, where $n=\frac{4^l-1}{3}$, with quotient matrix $(J-E)$. If $s$ is even then by Proposition~\ref{p:multiplyk} there is a perfect $4^l$-coloring in $D(m,n)$, where $2m+n=\frac{s(4^l-1)}{3}$, with quotient matrix $s(J-E)$ for any even $s$. 
     Let $s=2k+1$ is odd. By Proposition~\ref{p:2lham} there is a perfect $4^l$-coloring in $D(0,\frac{(k-1)(4^l-1)}{3})$ with quotient matrix $(k-1)(J-E)$.  
     By Proposition~\ref{p:multiplyk} there is a perfect $4^l$-coloring in $D(m,n)$, where $2m+n=\frac{2(k-1)(4^l-1)}{3}$, with quotient matrix $(2k-2)(J-E)$. By Corollary~\ref{c:3jham} there is a perfect $4^l$-coloring in $D(m,n)$, where $2m+n=4^l-1$, with quotient matrix $3(J-E)$ and there is a perfect $4^l$-coloring in $D(2^{2l-1},0)$ with quotient matrix $3J$. The required perfect coloring can be constructed by Theorem~\ref{t:diag}.   
\end{proof}

\begin{predln}\label{p:2lcol}
We have the following.
    \begin{enumerate}
        \item For any $l \ge 3$ and even $s$ there is perfect $2^{2l-1}$-coloring in $D(m,n)$, where $2m+n=\frac{s(2^{2l-1}-1)+a}{3}$, with quotient matrix $s(J-E)+aE$, where $a=0$ if ($s/2 =0 \mod 3)$, $a=2$ if $(s/2 =2 \mod 3)$ and $a=4$ if $(s/2=1 \mod 3)$.
        \item For any $l \ge 3$ and odd $s$ there is perfect $2^{2l-1}$-coloring in $D(m,n)$, where $2m+n=\frac{s(2^{2l-1}-1)+a}{3}$, with quotient matrix $s(J-E)+aE$, where $a=0$ if $((s-3)/2 =0 \mod 3)$, $a=2$ if $((s-3)/2 =2 \mod 3)$ and $a=4$ if $((s-3)/2=1 \mod 3)$.
        \item For any $l \ge 3$ and odd $s$ there is perfect $2^{2l-1}$-coloring in $D(m,0)$, where $2m=\frac{s(2^{2l-1}-1)+a}{3}$, with quotient matrix $s(J-E)+aE$, where $a=3$ if $((s-3)/2 =0 \mod 3)$, $a=5$ if $((s-3)/2 =2 \mod 3)$ and $a=7$ if $((s-3)/2=1 \mod 3)$.
        \item For $l=3$ and odd $s$ there is perfect $2^{2l-1}$-coloring in $D(m,n)$, where $2m=\frac{s(2^{2l-1}-1)+a}{3}$, with quotient matrix $s(J-E)+aE$, where $a=3$ if $((s-3)/2 =0 \mod 3)$, $a=5$ if $((s-3)/2 =2 \mod 3)$ and $a=7$ if $((s-3)/2=1 \mod 3)$.
    \end{enumerate}
\end{predln}
\begin{proof}
1) For $s=2$ the statement corresponds to Proposition~\ref{p:2lham}(i.4). 
Let $s=2r>2$.
By Proposition~\ref{p:2lham}(i.2) there is a perfect $2^{2l-1}$-coloring in $D(0,n')$, where $n'=\frac{r(2^{2l-1}-1)}{3}$, with quotient matrix $S=r(J-E)+aE$, where $a=0$ if $r =0 \mod 3$, $a=2$ if $r =2 \mod 3$ and $a=4$ if $r=1 \mod 3$. By Proposition~\ref{p:multiplyk} there is a perfect $2^{2l-1}$-coloring in $D(m,n)$, where $2m+n=\frac{2r(2^{2l-1}-1)}{3}$, with quotient matrix $2S$.

2-4) Let $s=2r+1$. Then the requirement perfect coloring can be obtained by Theorem~\ref{t:diag} using the perfect $2^{2l-1}$-colorings from item $1$ (for case $s-3$) and perfect $2^{2l-1}$-coloring from Corollary~\ref{c:3jham}.  
\end{proof}

\subsection{Constructions of (b,c)-colorings}

\begin{predln}\label{p:gcd}
    Let there is perfect $2^l$-coloring in $D(m,n)$ with quotient matrix $s(J-E)+aE$. Then for any $r \in \{1,2^l-1\}$ there is a perfect coloring in $D(m,n)$ with quotient matrix 
    $\displaystyle{\begin{pmatrix}
                s(r-1)+a & (2^l-r)s \\
                rs & (2^l--r-1)s
            \end{pmatrix}}$
\end{predln}
\begin{proof}
This coloring can be obtained if we unit first $r$ colors in one color and unit the other colors in the second color.  
\end{proof}

\begin{theoreman}
    Let $gcd(b,c)=1$ and $b+c=4^l$. Then there is a perfect coloring in $D(m,n)$, where $2m+n=\frac{4^l-1}{3}$, with the quotient matrix $\displaystyle{\begin{pmatrix}
                c-1 & b \\
                c & b-1
            \end{pmatrix}}$
\end{theoreman}
\begin{proof}
There is a $1$-perfect code in $D(m,n)$, where $2m+n=\frac{4^l-1}{3}$. The union of $c$ disjoint $1$-perfect codes corresponds to requirement perfect coloring.
\end{proof}

%\begin{theoreman}\label{t:bcex4l}
%    If $gcd(b,c)=1$ and $b+c=4^l$ then for any $n$ a pair $(b,c)$ is $n$-admissible. Moreover, $A_{n}(b,c) \le \frac{4^{l+1}-1-3n}{6}$.
%\end{theoreman}

\begin{theoreman}
Let $gcd(b,c) \ge 2$ and $\displaystyle{\frac{b+c}{gcd(b,c)}=4^l}$. Then 
        \begin{enumerate}
            \item There is a perfect coloring in $D(m,n)$, where $\displaystyle{2m+n=\frac{gcd(b,c)\cdot(4^l-1)}{3}}$, with quotient matrix $$\displaystyle{\begin{pmatrix}
                c-gcd(b,c) & b \\
                c & b-gcd(b,c)
            \end{pmatrix}}$$
            \item If $gcd(b,c)$ is odd then there is a perfect coloring in $D(m,0)$, where $\displaystyle{2m=\frac{gcd(b,c)(4^l-1)+3}{3}}$, with quotient matrix $$\displaystyle{\begin{pmatrix}
                c-gcd(b,c)+3 & b \\
                c & b-gcd(b,c)+3
            \end{pmatrix}}$$
    \end{enumerate}
\end{theoreman}
\begin{proof}
Let $b=gcd(b,c)b'$ and $c=gcd(b,c)c'$, where $b'$ and $c'$ are relatively prime.

1) This follows from Proposition~\ref{p:4lcol}(i.1) for case $s=gcd(b,c)$ and Proposition~\ref{p:gcd} for case $r=c'$.

2) This follows from Proposition~\ref{p:4lcol}(i.2) for case $s=gcd(b,c)$ and Proposition~\ref{p:gcd} for case $r=c'$.
%    Let $b=gcd(b)b'$ and $c=gcd(b,c)c'$, where $b'$ and $c'$ are relatively prime. By Proposition~\ref{p:4lcol} there is a $4^l$-perfect coloring $f$ in $D(m,n)$, where $2m+n=\frac{gcd(b,c)(4^l-1)}{3}$, with quotient matrix $gcd(b,c)(J-E)$. The requirement coloring can be obtained as union of first $c'$ colors of $f$. 
\end{proof}

%\begin{theoreman}
%     Let $gcd(b,c)=2$ and $\frac{b+c}{gcd(b,c)}=2^{2l-1}$. Then 
%     \begin{enumerate}
%         \item There is a perfect coloring in $D(m,n)$, where %$2m+n=\frac{2^{2l}-1}{3}$, with quotient matrix $\displaystyle{\begin{pmatrix}
%                c-1 & b \\
%                c & b-1
%            \end{pmatrix}}$  
%     \end{enumerate}        
%\end{theoreman}

\begin{theoreman}
    Let $gcd(b,c)=2$ and $\displaystyle{\frac{b+c}{gcd(b,c)}=2^{2l-1}}$ for some $l \ge 2$. Then 
    \begin{enumerate}
        \item There is a perfect coloring in $D(m,n)$, where $\displaystyle{2m+n=\frac{2^{2l}-1}{3}}$, with quotient matrix $\displaystyle{\begin{pmatrix}
                c-1 & b \\
                c & b-1
            \end{pmatrix}}$ 
        \item There is a perfect coloring in $D(m,0)$, where $\displaystyle{2m=\frac{2^{2l}+2}{3}}$, with quotient matrix $\displaystyle{\begin{pmatrix}
                c+2 & b \\
                c & b+2
            \end{pmatrix}}$ 
    \end{enumerate}
\end{theoreman}
\begin{proof}
Let $b=2b'$ and $c=2c'$, where $b'$ and $c'$ are relatively prime.

1) This follows from Proposition~\ref{p:2lham}(i.3) and Proposition~\ref{p:gcd} for case $r=c'$.

2) This follows from Proposition~\ref{p:2lham}(i.4) and Proposition~\ref{p:gcd} for case $r=c'$.
\end{proof}

\begin{theoreman}
    Let $gcd(b,c) > 2$ and $\displaystyle{\frac{b+c}{gcd(b,c)}=2^{2l-1}}$ for some $l \ge 2$. Then
        \begin{enumerate}
            \item If $gcd(b,c)$ is even then there is a perfect coloring in $D(m,n)$, where $\displaystyle{2m+n=\frac{gcd(b,c)(2^{2l-1}-1)+a}{3}}$, with quotient matrix $$\displaystyle{\begin{pmatrix}
                c-gcd(b,c)+a & b \\
                c & b-gcd(b,c)+a
            \end{pmatrix}},$$ where $a=0$ if $(gcd(b,c)/2 =0 \mod 3)$, $a=2$ if $(gcd(b,c)/2 =2 \mod 3)$, $a=4$ if $(gcd(b,c)/2 =1 \mod 3)$
            \item Let $gcd(b,c)$ is odd and if $l \ne 3$ then $m=0$. Then there is a perfect coloring in $D(m,n)$, where $\displaystyle{2m+n=\frac{gcd(b,c)(2^{2l-1}-1)+a}{3}}$, with quotient matrix $$\displaystyle{\begin{pmatrix}
                c-gcd(b,c)+a & b \\
                c & b-gcd(b,c)+a
            \end{pmatrix}},$$ where $a=0$ if $((gcd(b,c)-3)/2 =0 \mod 3)$, $a=2$ if $((gcd(b,c)-3)/2 =2 \mod 3)$, $a=4$ if $((gcd(b,c)-3)/2 =1 \mod 3)$
            \item If $gcd(b,c)$ is odd then there is a perfect coloring in $D(m,0)$, where $\displaystyle{2m=\frac{gcd(b,c)(2^{2l-1}-1)+a}{3}}$, with quotient matrix $$\displaystyle{\begin{pmatrix}
                c-gcd(b,c)+a & b \\
                c & b-gcd(b,c)+a
            \end{pmatrix}},$$ where $a=3$ if $((gcd(b,c)-3)/2 =0 \mod 3)$, $a=5$ if $((gcd(b,c)-3)/2 =2 \mod 3)$, $a=7$ if $((gcd(b,c)-3)/2 =1 \mod 3)$
             \item If $gcd(b,c)$ is odd and $l=3$ then there is a perfect coloring in $D(m,n)$, where $\displaystyle{2m+n=\frac{gcd(b,c)(2^{2l-1}-1)+a}{3}}$, with quotient matrix $$\displaystyle{\begin{pmatrix}
                c-gcd(b,c)+a & b \\
                c & b-gcd(b,c)+a
            \end{pmatrix}},$$ where $a=3$ if $((gcd(b,c)-3)/2 =0 \mod 3)$, $a=5$ if $((gcd(b,c)-3)/2 =2 \mod 3)$, $a=7$ if $((gcd(b,c)-3)/2 =1 \mod 3)$
        \end{enumerate}
\end{theoreman}
\begin{proof}
Let $b=gcd(b)b'$ and $c=gcd(b,c)c'$, where $b'$ and $c'$ are relatively prime. 
The statement follows from Proposition~\ref{p:2lcol} and Proposition~\ref{p:gcd} for case $r=c'$.
\end{proof}

\begin{predln}\label{p:bcpart}
Let $g^i_0,\ldots,g^i_{2^k-1}$, $i=1,\ldots,4$, are colorings (not necessary different) of $D(m,n)$ such that 
\begin{enumerate}
\item For any $i$ either $g^i_j$ is a perfect $(b_i,c_i)$-coloring
%with quotient matrix  
% $\displaystyle{\begin{pmatrix}
%a_i & b_i \\
%c_i & d_i 
%\end{pmatrix}}$
 For any $j$, either  $g^i_j \equiv 1$ for any $j$, or $g^i_j \equiv 2$ for any $j$;
 \item If $g^r_1$ and $g^s_1$ are not constant then $b_r+c_r=b_s+c_s=2^k$;
\item If $g^i_1$ is not constant  then any vertex of $D(m,n)$ has color $1$ in $c_i$ colorings among $g^i_1,\ldots,g^i_{2^k}$. 
\end{enumerate}
Then for any $k \ge 2$ there are perfect 
$(2^{k+2}-\gamma_1-\ldots-\gamma_4,\gamma_1+\ldots+\gamma_4)$-colorings $f_1,\ldots,f_{2^{k+2}}$ in $D(m+m',n+n')$, where $2m'+n'=2^k$ (if $k=3$ then $m'=0$), 
such that any vertex of $D(m+m',n+n')$ has color $1$ in $\gamma_1+\ldots+\gamma_4$ colorings among $f_1,\ldots,f_{2^{k+2}}$, where $\gamma_i=c_i$ if $g^i_1$ is not constant, $\gamma_i=2^k$ if $g^i_1 \equiv 1$, $\gamma_i=0$ if $g^i_1 \equiv 2$. 
\end{predln}
\begin{proof}
Denote numbers $1,\ldots,2^{k+2}$ as $(i,j)$, where $(i,j)=i2^k+j$,  $i \in \{0,\ldots,3\}$, $j \in \{0,\ldots,2^k-1\}$. 
By Corollary~\ref{c:multipart} 
there is $k$-multipart perfect coloring $h$ in $D(m',n')$ in colors $(0,0),\ldots,(3,2^k-1)$.
%graph $D(m',n')$, where $2m'+n'=2^k$, can be partition %into $2^{k+2}$ disjoint  codes $C^(0,0),\ldots,C^(2^k-%1,2^k-1)$, such that for any 
%$i \ne j$ and $r,s$ any vertex from $C^{(i,r)}$ adjacent %to exactly $1$ vertex form $C^{(j,s)}$ and for any %$i$,$r$ and $s$ any vertex from $C^{(i,r)}$ has no %neighbours in $C^{(i,s)}$.
Represent a vertices of $D(m'+m,n'+n)$ as $(x,y)$, where $x \in \VV(D(m',n'))$, $y \in \VV(D(m,n))$. 

Define function $f^{(i,j)}$ in the following way:
$$f^{(i,j)}(x,y)=g^{r+i}_{s+j}(y) \text{ if } h(x)=(r,s)$$
where sum $(r+i)$ is modulo $4$ and $(s+j)$ is modulo $2^k$.

Consider coloring $f^{(i,j)}$.
Let $(x,y)$ be an arbitrary vertex of $D(m'+m,n'+n)$ such that $h(x)=(r,s)$. Denote by $N_1(x,y)$ the set of vertices $(x,y')$, where $y'$ adjacent to $y$ in $D(m,n)$. Denote by $N^i_2(x,y)$ the set of vertices $(x',y)$, where $x'$is a vertex adjacent to $x$ in $D(m',n')$ such that $h(x')=(i,l)$ for some $l$. The neighbourhood $N(x,y)$ of $(x,y)$ equals $N(x,y)=N_1(x,y) \cup N^0_2(x,y) \cup \ldots \cup N^3_2(x,y)$. 

Let $(x,y)$ has color $1$ in coloring $f^{(i,j)}$.
The number of vertices in $N_1(x,y)$ that have color $2$ equals the number of vertices $y'$ adjacent to $y$ that have color $2$ in coloring $g^{r+i}_{s+j}$. If this coloring is not constant then this number equals $b_{r+i}=2^k-c_{r+i}=2^k-\gamma_{r+i}$. If $g^{r+i}_{s+j}\equiv 1$ (case $g^{r+i}_{s+j}\equiv 2$ is not possible) then this number equals $0$ (in this case, this equals $2^k-\gamma_{r+i}$). 
The number of $(x',y)$ from $N^t_2(x,y)$, $t \ne i+r$, that has color $2$ equals the number of colorings among $g^t_0,\ldots,g^t_{2^k-1}$ such that $y$ in this coloring has color $2$. By condition this number equals $2^k-\gamma_t$.
If $t=i+r$ the set $N^t_2(x,y)$ is empty.
Therefore, the number of neighbours of $(x,y)$ that have color $2$ equals $2^{k+2}-\gamma_1-\ldots-\gamma_4$.

Let $(x,y)$ has color $2$ in coloring $f^{(i,j)}$.
The number of vertices in $N_1(x,y)$ that have color $1$ equals the number of vertices $y'$ adjacent to $y$ that have color $1$ in coloring $g^{r+i}_{s+j}$. If this coloring is not constant then this number equals $c_{r+i}=\gamma_{r+i}$. If $g^{r+i}_{s+j}\equiv 2$ (case $g^{r+i}_{s+j}\equiv 1$ is not possible) then this number equals $0$ (in this case this equals $\gamma_{r+i})$. 
The number of $(x',y)$ from $N^t_2(x,y)$, $t \ne i+r$, that has color $1$ equals the number of colorings among $g^t_0,\ldots,g^t_{2^k-1}$ such that $y$ in this coloring has color $1$. By condition this number equals $\gamma_t$.
If $t=i+r$ the set $N^t_2(x,y)$ is empty.
Therefore, the number of neighbours of $(x,y)$ that have color $1$ equals $\gamma_1+\ldots+\gamma_4$.

Therefore, for any $i,j$ the coloring $f^{(i,j)}$ is perfect with requirement parameters.

It remains to prove that arbitrary vertex $(x,y)$ such that $h(x)=(r,s)$ has color $1$ in $\gamma_1+\ldots+\gamma_4$ colorings among $f^{(0,0)},\ldots,f^{(0,2^k-1)}$.
The vertex has color $1$ in $f^{(i,j)}$ if and only if $y$ has color $1$ in coloring $g^{r+i}_{j+s}$. Since indexes of $g$ runs through $(0,0),\ldots,(0,2^k-1)$ the number of $(i,j)$ such that $f^{(i,j)}(x,y)=1$ equals the number of $(l,t)$ such that $y$ has color $1$ in $g^l_t$. This number equals $\gamma_1+\ldots+\gamma_4$. 

\end{proof}

\begin{figure}
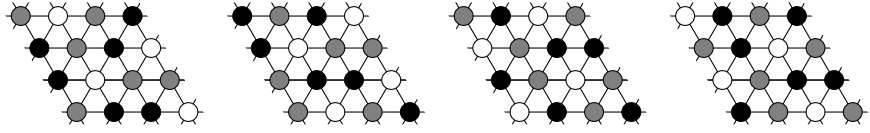

\centering
 \raisebox{-12mm}{\sh{xz nz xz zz}{zz xz zz nz}{zz nz xz xz}{xz zz zz nz}}\
 \raisebox{-12mm}{\sh{zz xz zz nz}{zz nz xz xz}{xz zz zz nz}{xz nz xz zz}}\
 \raisebox{-12mm}{\sh{xz zz nz xz}{nz xz zz zz}{zz xz nz xz}{nz zz xz zz}}\ 
 \raisebox{-12mm}{\sh{nz zz xz zz}{xz zz nz xz}{nz xz zz zz}{zz xz nz xz}}
\caption{Illustration to Lemma~\ref{l:bcind}, item~1.}
\label{f:1533part}

\end{figure}

\begin{figure}
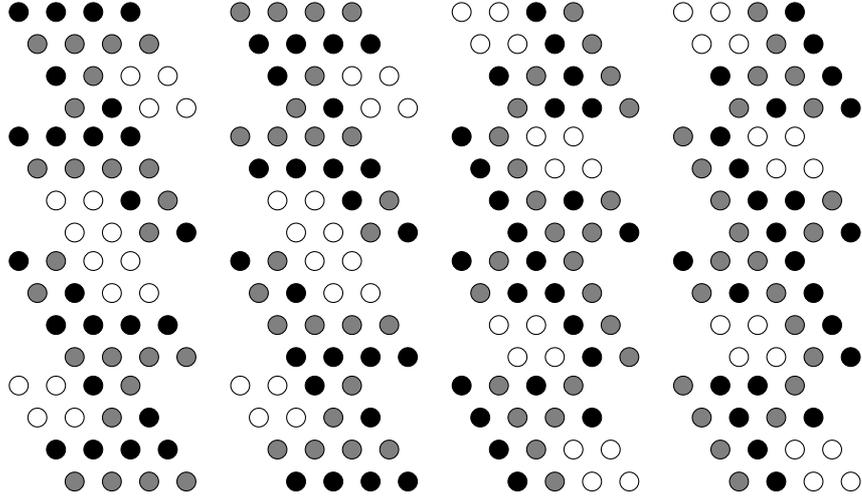

\centering
 \raisebox{-12mm}{\hh{zz zz zz zz}{xz xz xz xz}{zz xz nz nz}{xz zz nz nz}}\ 
 \raisebox{-12mm}{\hh{xz xz xz xz}{zz zz zz zz}{zz xz nz nz}{xz zz nz nz}}\
\raisebox{-12mm}{\hh{nz nz zz xz}{nz nz zz xz}{zz xz zz xz}{xz zz zz xz}}\ 
 \raisebox{-12mm}{\hh{nz nz xz zz}{nz nz xz zz}{zz xz xz zz}{xz zz xz zz}}

 \raisebox{-12mm}{\hh{zz zz zz zz}{xz xz xz xz}{nz nz zz xz}{nz nz xz zz}}\ 
 \raisebox{-12mm}{\hh{xz xz xz xz}{zz zz zz zz}{nz nz zz xz}{nz nz xz zz}}\
\raisebox{-12mm}{\hh{zz xz nz nz}{zz xz nz nz}{zz xz zz xz}{zz xz xz zz}}\ 
 \raisebox{-12mm}{\hh{xz zz nz nz}{xz zz nz nz}{xz zz zz xz}{xz zz xz zz}}

 \raisebox{-12mm}{\hh{zz xz nz nz}{xz zz nz nz}{zz zz zz zz}{xz xz xz xz}}\ 
 \raisebox{-12mm}{\hh{zz xz nz nz}{xz zz nz nz}{xz xz xz xz}{zz zz zz zz}}\
\raisebox{-12mm}{\hh{zz xz zz xz}{xz zz zz xz}{nz nz zz xz}{nz nz zz xz}}\ 
 \raisebox{-12mm}{\hh{zz xz xz zz}{xz zz xz zz}{nz nz xz zz}{nz nz xz zz}}

 \raisebox{-12mm}{\hh{nz nz zz xz}{nz nz xz zz}{zz zz zz zz}{xz xz xz xz}}\ 
 \raisebox{-12mm}{\hh{nz nz zz xz}{nz nz xz zz}{xz xz xz xz}{zz zz zz zz}}\
\raisebox{-12mm}{\hh{zz xz zz xz}{zz xz xz zz}{zz xz nz nz}{zz xz nz nz}}\ 
 \raisebox{-12mm}{\hh{xz zz zz xz}{xz zz xz zz}{xz zz nz nz}{xz zz nz nz}}
\caption{Illustration to Lemma~\ref{l:bcind}, item~2.}
\label{f:4536part}
\end{figure}

\begin{lemman}\label{l:bcind}
We have the following.
\begin{enumerate}
    \item For any $c\in \{3,\ldots,6\}$ there are $8$ perfect $(8-c,c)$-colorings $g_1,\ldots,g_8$ in $D(1,0)$ such that any vertex of $D(1,0)$ has color $1$ in $c$ colorings among $g_1,\ldots,g_8$.
    \item For any $c \in \{3,\ldots,6\}$ there are $8$ perfect $(8-c,c)$-colorings $h_1,\ldots,h_8$ in $H(3,4)$ such that any vertex of $H(3,4)$ has color $1$ in $c$ colorings among $h_1,\ldots,h_8$. 
\end{enumerate}
\end{lemman}
\begin{proof}
    1) On Figure~\ref{f:1533part} is shown $8$ completely regular codes with intersection array $(5,3)$ (black and grey color),  denote these codes as $A_1,\ldots,A_8$, and $4$ completely regular codes with intersection array $(6,2)$ (white color), denote these codes $B_1,\ldots,B_4$. 
    If $c=2$ we can take $2$ times each code among $B_1,\ldots,B_4$.
    If $c=3$ we can take $1$ times each code among $C_1,\ldots,C_8$.
    If $c=4$ we can take $4$ times each code among $B_1 \cup B_2, B_3 \cup B_4$.
    If $c \in \{5,6\}$ we can change colors in case $c \in \{2,3\}$.

    2) On Figure~\ref{f:4536part} is shown $8$ completely regular codes with intersection array $(5,3)$ (black and grey color) and $4$ completely regular codes with intersection array $(6,2)$ (white color). In this case constructions of requirement codes are analogous to construction in previous item. 
\end{proof}

\begin{theoreman}\label{t:bc2l}
    Let $gcd(b,c)=1$, $b,c \ne 1$ and $b+c=2^{2l-1}$, $l \ge 3$.
    Then
    \begin{enumerate}
        \item There is a perfect coloring in $D(m,n)$, where $2m+n=\frac{b+c-2}{3}$, $n \ge 8$, $m \ge 1$ with quotient matrix  $\displaystyle{\begin{pmatrix}
c-2 & b \\
c & b-2 
\end{pmatrix}}$.
\item There is a perfect coloring in $H(n,4)$, where $n=\frac{b+c+1}{3}$, with quotient matrix  $\displaystyle{\begin{pmatrix}
c+1 & b \\
c & b+1 
\end{pmatrix}}$
    \end{enumerate}
\end{theoreman}
\begin{proof}
   1) We will prove by induction on $l$ that for any $c \in \{2,\ldots,2^{2l-1}-2\}$ there are $2^{2l-1}$ perfect $(2^{2l-1}-c,c)$-colorings in $D(m,n)$, where $2m+n=\frac{b+c-2}{3}$, such that any vertex of $D(m,n)$ has color $1$ in $c$ colorings.  Case $l=2$ considered in Lemma~\ref{l:bcind} (i.1). Suppose the statement true for some $l$. Let us prove the statement for $l+1$. 
   
   Note, that $c$ can be represented as $\gamma_1+\ldots+\gamma_4$, where $\gamma_i \in \{0,1,\ldots,2^k-1\}$, $\gamma_i \ne 1, 2^{k}-1$. Indeed, let $c=x \mod 2^k$. If $x$ not equals $1$ or $2^k-1$ modulo $2^k$ then it can be represented as $\gamma_1=x$ and $\gamma_i$, $i \ne 1$, equals $0$ or $2^k$. If $x=1$ then $\gamma_1=2^k-2$, $\gamma_2=3$ and $\gamma_i$, where $i=3,4$, equals $0$ or $2^k$. If $x=2^k-1$ then $\gamma_1=2^k-4$, $\gamma_2=3$ and $\gamma_i$, where $i=3,4$, equals $0$ or $2^k$. Note that on the step ($l=2 \to l=3$) we use $3$-multipart perfect coloring that is known for $D(m,n)$ only in case $m=0$ and $n=8$. So, we have $n \ge 8$.

   By induction hypothesis if $c \in \{2,\ldots.2^k-2\}$  then there are $2^{2l-1}$ perfect $(2^{2l-1}-c,c)$-colorings in $D(m,n)$, where $2m+n=\frac{b+c-2}{3}$, such that any vertex of $D(m,n)$ has color $1$ in $c$ colorings. Hence, the statement for $l+1$ follows from Proposition~\ref{p:bcpart}.

2) The item~2 can be proved analogously from Lemma~\ref{l:bcind} (i.2) and Proposition~\ref{p:bcpart}. 
\end{proof}

\begin{predln}
There is a perfect $(b,b)$-coloring in $D(m,n)$ if and only if $b$ is even and $2b=4i$ for some $i \in \{1,\ldots,2m+n\}$.   
\end{predln}
\begin{proof}
    The necessity follows from Proposition~\ref{p:neceigen}. 
    
    There is a perfect $(4,4)$-coloring (the set $\{00,01,02,03,20,21,22,23\}$ forms the first color) and perfect $(2,2)$-coloring (the set $\{00,01,02,03,10,11,12,13\}$ forms the first color) in $D(1,0)$. Also there is a perfect $(2,2)$-coloring in $D(0,1)$ (the set $\{0,1\}$ forms the first color). 

    Let us prove the statemeте by induction on $k$ that equals diameter $2m+n$. For $k=1,2$ the statement is true. Suppose that the statement is true for all $l \le k$ for some $k$. Consider graph $D(m,n)$ of diameter $2m+n=k+1$. If $n>0$ then by induction hypothesis there is a perfect $(b',b')$-coloring, where $b' \in \{1,\ldots,\ 4m+2n-2\}$. If $b \in \{1,\ldots,4m+2n-2\}$ there is a perfect $(b,b)$-coloring in $D(m,n)$ by Proposition~\ref{p:admiss}. If $b=4m+2n$ there is a perfect $(b,b)$-coloring in $D(m,n)$ by Theorem~\ref{t:diag} using a perfect $(2,2)$-coloring. If $n=0$ the requirement coloring can be obtained analogously by using Proposition~\ref{p:admiss}, Theorem~\ref{t:diag} and perfect colorings in $D(1,0)$ mentioned previously.  
\end{proof}

Summarizing the necessary conditions from Proposition~\ref{p:neceigen} and the results of the current section we have the following. 

\begin{theoreman}\label{t:admissible}
We have the following
    \begin{enumerate}
      \item A pair $(b,1)$ is $\infty$-admissible if and only if $b=4^l-1$ for some $l$. 
      \item A pair $(b,c)$, where $b \ne 1$ and $c \ne 1$, is $\infty$-admissible if and only if $\displaystyle{\frac{b+c}{gcd(b,c)}=2^l}$ for some $l$.
      \item A pair $(b,1)$ is $1$-admissible if and only if $b=4^l-1$ for some $l$.
      \item If $gcd(b,c) \ge 2$ and $\displaystyle{\frac{b+c}{gcd(b,c)}=2^l}$ for some $l$  then a pair $(b,c)$ is $0$-admissible.
      \item If $gcd(b,c)=1$ and $\displaystyle{\frac{b+c}{gcd(b,c)}=4^l}$ for some $l$  then a pair $(b,c)$ is $1$-admissible.
      \item If $gcd(b,c)=1$ and $\displaystyle{\frac{b+c}{gcd(b,c)}=2^{2l-1}}$ for some $l$  then a pair $(b,c)$ is $8$-admissible.
    \end{enumerate}
\end{theoreman}

%\begin{theoreman}\label{t:bcex2l}
%\begin{enumerate}
%    \item If $gcd(b,c) \ge 2$ then a pair $(b,c)$ is a complete admissible
%    \item If $gcd(b,c)=1$ and $b+c=2^{2l-1}$ then $(b,c)$ is a complete %admissible intersection array
%     \item If $gcd(b,c)=1$ and $b+c=4^l$ then $[b,c]$ is a $h$-admissible %intersection array.   
%\end{enumerate}
%\end{theoreman}

%\begin{predln}\label{p:hamcover}
%\begin{enumerate}
%\item There is such $a,m,n,N$ that the graph $D(m,n)$ cover graph $lH(N,2^K)+aI$. 
%\item The graph $D(m,n)$ cover graph $H(N,4^k)$, where $2m+n=N\frac{4^l-1}{3}$, $n \ge N$.
%\end{enumerate}
%\end{predln}

%\begin{coroll}
%\begin{enumerate}
%    \item If there is a completely regular code with intersection array $[\beta_0,\ldots, \beta_{\rho-1};\gamma_1,\ldots,\gamma_{\rho}]$ in $H(N,2^k)$ for some $N$ and $k$ then for any $l \ge 2$ the intersection array $[l\beta_0,\ldots, l\beta_{\rho-1};l\gamma_1,\ldots,l\gamma_{\rho}]$ is complete admissible.
%    \item If there is a perfect coloring in $H(N,2^k)$ with quotient matrix then there is a perfect coloring in $D(m,n)$, where  $2m+n=N\frac{4^l-1}{3}$, $n \ge N$, with quotient matrix $S$
%\end{enumerate}
    
%\end{coroll}

\section{Multifold perfect codes}
In \cite{Kro:multperf} it was obtained the characterization of parameters of multifold $1$-perfect codes in Hamming graphs on the prime power alphabet.

\begin{theoreman}\cite{Kro:multperf}
    Let $q=p^t$ for a prime $p$. A $\mu$-fold $1$-perfect code in $H(n,q)$ 
    exists if and only if 
\begin{enumerate}
    \item[(i)]  $\mu<n(q-1)+1$,
    \item[(ii)] $\displaystyle{\frac{\mu q^n}{n(q-1)+1}}$ is integer,
    \item[(iii)] and $n \equiv 1 \mod q$.
\end{enumerate}
\end{theoreman}

In current paper, we obtain the same result for the Doob graphs.

\begin{predln}\label{p:tiling}
Let there is a perfect $2^l$-coloring $g$ in $D(m,n)$ with quotient matrix
$k(J-E)+aE$. Let there are $2^l$ perfect $2$-colorings in $D(m',n')$ with quotient matrix 
$\displaystyle{\begin{pmatrix}
a' & b' \\
c' & d' 
\end{pmatrix}}$
such that any vertex of $D(m',n')$ has color $1$ in exactly $r$ colorings.
Then there are $2^l$ perfect coloring in $D(m+m',n+n')$ with quotient matrix 
$$\displaystyle{\begin{pmatrix}
a'+a+(r-1)k & b'+(2^l-r)k \\
c'+rk & d'+(2^l-r-1)k+a 
\end{pmatrix}}$$
such that any vertex of $D(m+m',n+n')$ has color $1$ in exactly $k$ colorings
\end{predln}
\begin{proof}
Denote the vertices of $D(m+m',n+n')$ as $(x,y)$, where $x \in \VV{D(m,n)}$ and $y \in \VV{D(m',n')}$. 
Denote by $N_1(x,y)$ the set of vertices $(x,y')$, where $y'$ adjacent to $y$.
Denote by $N_2(x,y)$ the set of vertices $(x',y)$, where $x'$ adjacent to $x$.
The neighbourhood of $(x,y)$ equals $N(x,y)=N_1(x,y) \cup N_2(x,y)$.
Denote colorings of $D(m'.n')$ as $f_1,\ldots,f_{2^l}$.

Define coloring $f^i$, where $i \in \{1,\ldots,2^l\}$, of $D(m+m',n+n')$ in the following way:
$$f^i(x,y)=f_{g(x)+i-1}(y),$$
where sum $g(x)+i-1$ is modulo $2^l$.

Let us prove that for arbitrary $i$ coloring $f^i$ is perfect with requirement quotient matrix. 

%Let $(x,y)$ be arbitrary vertex of $D(m+m',n+n')$, where $g(x)=j$.

Let $(x,y)$ has color $1$.
The number of vertices in $N_1(x,y)$ that have color $1$ equals $a'$.
The number of vertices in $N_2(x,y)$ that have color $1$ equals the number of such vertices $z$ that it is adjacent to $x$ in $D(m,n)$ and $f_{g(z)+i-1}(y)=1$. This equals sum $a$ (when $g(z)=g(x)$) and $(r-1)k$ ($k$ vertices for each $z$ where  $g(z) \ne g(x)$ and $f_{g(z)+i-1}(y)=1$).

The case when $(x,y)$ has color $2$ is analogous.
The number of vertices in $N_1(x,y)$ that have color $1$ equals $c'$.
The number of vertices in $N_2(x,y)$ that have color $1$ equals $rk$ ($k$ vertices for each $z$ where $f_{g(z)+i-1}(y)=1$).

It remains to prove that $(x,y)$ has color $1$ in $k$ colorings among $f^1,\ldots,f^{2^l}$.
The vertex $(x,y)$ has color $1$ in $f^i(x,y)$ if and only if $f_{g(x)+i-1}(y)=1$. Since index $(g(x)+i-1)$ runs through all elements from $\{1,\ldots,2^l\}$ the number of $i$ such that $f^i(x,y)=1$ equals the number of $j$ such that $f_j(y)=1$, which equals $r$.  
\end{proof}

\begin{theoreman}\label{t:mlrtipc}
There is a $c$-fold perfect code in $D(m,n)$ if and only if 
\begin{enumerate}
\item $2m+n=1 \mod 4$
\item $c=\alpha l$, where $\alpha$ is such that $6m+3n+1=\alpha 2^k$, $\alpha$ is odd and $l \in \{1,\ldots,2^k\}$
\end{enumerate}
\end{theoreman} 
\begin{proof}
The necessary condition follows from Lloyd conditions (the graph should have  eigenvalue $-1$) and the fact that cardinality of such code equals $\displaystyle{\frac{c \cdot 4^{2m+n}}{(6m+3n+1)}}$. 
 
By necessary conditions multifold $1$-perfect code can exist only if $2m+n=4k+1$.  
So we prove the statement by induction on $k$. 
Moreover, we prove that for any $k$ there is a partition of $D(m,n)$, $2m+n=4k+1$, into $2^l$ disjoint $\alpha$-fold $1$-perfect codes, where $3(2m+n)+1=\alpha  2^l$, $\alpha$ is odd.

For $k=0$ there is a partition of $D(0,1)$ into $4$ disjoint $1$-perfect codes.
For $k=1$ there is a partition of $D(2,1)$, $D(1,3)$ and $D(0,5)$ into $16$ disjoint $1$-perfect codes.

Suppose we prove the induction hypothesis for all $k'<k$. Let us prove it for $k$.
Consider $D(m,n)$ of diameter $D=2m+n=4k+1$. 
Represent $D$ as $2^r+D'$, where $D'<2^r$. We have $3D+1=3 \cdot 2^r+3D'+1$. 
Represent $D(m,n)$ as direct product $D(m'',n'')$, where $2m''+n''=2^r$, and $D(m',n')$, where $2m'+n'=D'$. 
Let $3D'+1=\beta \cdot 2^s$, where $\beta$ is odd.

Consider some cases.

1) The case $(s<r)$. We have $3D+1=(3 \cdot 2^{r-s}+\beta)2^s$. So we should prove that there is a partition of $D(m,n)$ into $2^s$ disjoint $(3 \cdot 2^{r-s}+\beta)$-fold $1$-perfect codes. By Corollary~\ref{c:3jham} there is a perfect $2^r$-coloring in $D(m'',n'')$ with quotient matrix $3J$. If we partition colors into $2^{s}$ blocks that consist of $2^{r-s}$ colors and unit colors in each block we obtain perfect $2^s$-coloring with quotient matrix $3 \cdot 2^{r-s}(J_{2^s \times 2^s})$. By induction hypothesis there are $2^s$ disjoint $\beta$-fold $1$-perfect codes in $D(m',n')$.
By Proposition\ref{p:tiling} there are $2^s$ disjoint $(\beta+3 \cdot 2^{r-s})$-fold 
$1$-perfect codes.

2) The case $(s>r)$. Since $D'<2^r$ we have $3D'+1<2^{r+2}$, hence, $s=r+1$ and $\beta<2$. Note that in this case $3D+1=2^{r+1}$ and, hence, $r$ should be odd.
In this case, there is a partition of $D(m',n')$ into $2^{k+1}$ disjoint $1$-perfect codes.
So, we should prove that there is a partition of $D(m,n)$ into $2^r$ disjoint $5$-fold $1$-perfect codes. By Corollary~\ref{c:3jham} there is a perfect $2^r$-coloring in $D(m'',n'')$ with quotient matrix $3J$. 
In $D(m',n')$ there is a partition into $2^r$ disjoint $2$-fold $1$-perfect codes. Hence, by Proposition\ref{p:tiling} there are $2^r$ disjoint $5$-fold 
$1$-perfect codes.

3) The case $(r=s)$. Since $D'<2^r$ we have $3D'+1<2^{r+2}$. Hence, we have $\beta=1$ or $\beta=3$. If $\beta=3$ we have a contradiction with $3D'+1=3 \cdot 2^r$. So, $\beta=1$ and $3D+1=3\cdot2^r+2^{r}=2^{r+2}$. Hence, $D=\frac{2^{r+2}-1}{3}$. If $r$ is odd this is not integer. If $r$ is even then in this case there is a partition of $D(m,n)$ into $2^{r+2}$ disjoint $1$-perfect codes.

Finally, note that we can choose any $m''$ and $n''$ such that $2m''+n''=2^r$ and any $m'$ and $n'$ (by induction hypothesis) such that $2m'+n'=D'$, so, the statement is true for any $m$ and $n$ such that $2m+n=D$.
\end{proof}

\bibliographystyle{unsrt}
\bibliography{k}

\end{document}